\documentclass[11pt]{article}
\usepackage{amsfonts,amsmath,amssymb,amsthm, amscd, subcaption}
\usepackage{graphicx}
\usepackage[stable]{footmisc}
\numberwithin{equation}{section}
\usepackage{xfrac}
\usepackage{tikz}
\usepackage{tikz-cd}
\usetikzlibrary{decorations.markings}
\usetikzlibrary{arrows}
\usetikzlibrary{calc}

\newtheorem{theorem}{Theorem}[section]
\newtheorem{prop}[theorem]{Proposition}
\newtheorem{lemma}[theorem]{Lemma}
\newtheorem{cor}[theorem]{Corollary}

\theoremstyle{definition}
\newtheorem{definition}[theorem]{Definition}

\newtheorem{remark}[theorem]{Remark}

\usepackage{float, mathtools}

\newcommand{\Mr}{M_A^{\textup{red}}}

\def\({\langle \hskip -.1cm \langle} 
\def\){\rangle \hskip -.1cm \rangle} 
\def\<{{\langle}}
\def\>{{\rangle}}

\def\a{{\alpha}}

\def\g{{\gamma}}
\def\i{{\iota}}

\def\d{{\delta}}
\def\Z{\mathbb Z}
\def\N{\mathbb N}

\def\S{\mathbb S}

\def\P{{\cal P}}
\def\Z{\mathbb Z}
\def\a{\alpha}

\def\k{{\kappa}}

\def\x{\chi}

\def\la{\lambda}
\def\M{{\cal M}}

\def\s{\sigma}
\def\e{\epsilon}

\def\Aut{{\textup{Aut}}}

\def\ni{\noindent} 
\begin{document}

\title{Peripheral structures of core groups} 

\author{Daniel S. Silver
\and
Lorenzo Traldi}

\date{}

\maketitle %{\setlength{\linewidth}{2in}

%%%%%%%%%%%%%%%%%%%%%%%%%%%%%% 

\begin{abstract} 
The core group is an invariant of unoriented virtual links. We introduce a peripheral structure for the core group, in which the longitudes are sensitive to orientations. We show that the combination of the core group and its peripheral structure is equivalent, as a link invariant, to the combination of the $\pi$-orbifold group and its peripheral structure. Examples show that the peripheral structure of the core group can be used to verify noninvertibility of some knots and links. 
\end{abstract}

\begin{center} Keywords: \textit {core group; link; longitude; meridian; peripheral system; $\pi$-orbifold group}

 \end{center}

\section{Introduction} \label{Intro} 

We follow the usual conventions regarding virtual link diagrams. To wit: a virtual link diagram $D$ is built from a finite number $\mu$ of piecewise smooth, oriented closed curves in the plane. The curves must be in general position, i.e., the only (self-)intersections are transverse double points, called crossings. At a classical crossing, one incident curve segment is designated the underpasser, and small portions of it are removed on each side of the crossing. At a virtual crossing, we draw a small circle around the crossing and formally view the crossing as ``not really being there.'' That is, we pretend the incident curve segments do not intersect. Removing the small portions of underpassing arcs at classical crossings has the effect of cutting the original curves into arcs. The set of arcs of $D$ is denoted $A(D)$, and the set of classical crossings of $D$ is denoted $C(D)$. If $c$ is a classical crossing we use $a(c)$ to denote the overpassing arc at $c$, and $b_1(c),b_2(c)$ to denote the underpassing arc(s) at $c$, in no particular order. A (virtual) link is an equivalence class of virtual link diagrams under detour moves and Reidemeister moves. If a virtual link diagram has no virtual crossing, the diagram and the associated link are classical. 

A.J. Kelly \cite{kelly} and M. Wada \cite{W} introduced the following notion into classical knot theory over thirty years ago. 

\begin{definition}
\label{arccore}
 The \emph{core group} of a virtual link diagram $D$ is the group $AC(D)$ given by the presentation 
\[\langle \{g_a \mid a \in A(D)\};\{r_c \mid c \in C(D)\} \rangle\text{,}\]
where the relator asssociated to a  classical crossing $c \in C(D)$ is $r_c=g_{a(c)}g^{-1}_{b_1(c)}g_{a(c)}g^{-1}_{b_2(c)}$.
\end{definition}

In this paper we abuse notation by using the symbols $g_a$ to represent elements of $AC(D)$, rather than adopting explicit notation for the epimorphism mapping a free group onto $AC(D)$ that arises from Definition \ref{arccore}.

In earlier work with S. G. Williams \cite{STW1}, we discussed the connection between $AC(D)$ and a couple of related groups defined using the regions of a diagram. The purpose of the present paper is to introduce a peripheral structure for $AC(D)$. This peripheral structure is tied to a particular family of automorphisms of $AC(D)$.

\begin{definition} \label{immap}
If $D$ is a virtual link diagram then for each $a \in A(D)$ there is an automorphism $\i_{g_a}$ of $AC(D)$ given by the formula
\[
\i_{g_a}(g_b) = g_a g_b^{-1} g_a \thickspace \forall b \in A(D).
\]
The subgroup of the automorphism group $\Aut(AC(D))$ generated by these automorphisms is denoted $\mathcal{I}(D)$.
\end{definition}

We say $\i_{g_a}$ is the \emph{involutory meridional automorphism} associated with $g_a$. The name is clumsy but accurate: it reflects the facts that $\i_{g_a}$ is determined by the meridional element $g_a$, and is an involution; that is, $\i_{g_a} \circ \i_{g_a}$ is the identity map of $AC(D)$. For convenience we often refer to involutory meridional automorphisms as $\i_{g_a}$ maps, rather than using the full name.

\begin{definition} \label{meridian}
Let $D$ be a diagram of $L=K_1 \cup \cdots \cup  K_ \mu$. Then the \emph{meridians} in $AC(D)$ are defined recursively as follows.
\begin{enumerate}
    \item If $b \in A(D)$ is an arc of $K_i$ in $D$ then $g_b$ is a meridian of $K_i$ in $AC(D)$.
    \item If $m$ is a meridian of $K_i$ in $AC(D)$ and $a \in A(D)$ then $\i_{g_a}(m)$ is also a meridian of $K_i$ in $AC(D)$.
\end{enumerate}
The set of meridians in $AC(D)$ is denoted $\M(D)$, and if $1 \leq i \leq \mu$ then $\M_i(D)$ is the set of meridians associated with $K_i$.
\end{definition}

We should emphasize that we use the word ``meridian'' to refer to an element of $AC(D)$, not a curve in any topological space.

Let $\kappa_D:\M(D) \to \{1, \dots, \mu \}$ be the function that assigns each meridian the index of its associated component. Also, let $A_\mu = \Z \oplus (\Z_2)^{\mu-1}$. Then it is not hard to verify that there is a well-defined homomorphism $\phi_D:AC(D) \to A_\mu$ defined as follows. If $m \in \M_1(D)$ then $\phi_D(m) = (1,0, \dots, 0)$. If $i>1$ and $m \in \M_i(D)$ then $\phi_D(m) = (1,0, \dots, 0, 1 ,0, \dots, 0)$, with the second 1 in the $i$th coordinate. Notice that permuting the indices of $K_1, \dots, K_\mu$ will change both $\k_D$ and $\phi_D$.

Orientations are not used in defining the group $AC(D)$, or the meridians. We do use orientations to define longitudes, though.  Like the meridians, the longitudes in $AC(D)$ are defined in two steps. 

The first step is this. Suppose we are given an arc $b \in A(D)$, with $\k_D(b)=i$. We index the arcs of $K_i$ as $b=b_{i0},b_{i1},\dots, b_{i(k_i-1)},b_{ik_i}=b$, in order according to the orientation of $K_i$. For $0 \leq j \leq k_i$ let $c_{ij}$ be the crossing of $D$ separating $b_{ij}$ from $b_{i(j+1)}$, and let $a_{ij}$ be the overpassing arc at the crossing $c_{ij}$. Then we define the longitude $\la(b)$ as follows.
\[
\la(b) = \begin{cases}
g_{a_{i0}}g_{a_{i1}}^{-1} \dots g_{a_{i(k_i-2)}} g_{a_{i(k_i-1)}}^{-1} \text{,} &\text{  if }k_i\text{ is even} \\
g_{a_{i0}}g_{a_{i1}}^{-1} \dots g_{a_{i(k_i-2)}}^{-1} g_{a_{i(k_i-1)}} g_b^{-1}\text{,} & \text{ if }k_i\text{ is odd}
\end{cases}
\]

As with the term ``meridian,'' we use the term ``longitude'' to refer to a group element, not a curve. Observe that every longitude $\la(b)$ is represented by a word with alternating exponents $+1,-1,+1, \dots , -1$. 

Observe also that the definition gives the following formulas for the longitude associated with the arc $b_{i1}$ that meets $b=b_{i0}$ at the crossing $c_{i0}$.

\[
\la(b_{i1}) = \begin{cases}
g_{a_{i1}}g_{a_{i2}}^{-1} \dots g_{a_{i(k_i-1)}} g_{a_{i0}}^{-1} \text{,} &\text{  if }k_i\text{ is even} \\
g_{a_{i1}}g_{a_{i2}}^{-1} \dots g_{a_{i(k_i-1)}}^{-1} g_{a_{i0}} g_{b_{i1}}^{-1}\text{,} & \text{ if }k_i\text{ is odd}
\end{cases}
\]
As $g_{a_{i0}} g_{b_{i1}}^{-1}=g_b g_{a_{i0}}^{-1} $, we see that $\la(b_{i1}) = g_{a_{i0}}\overline{\la(b)} g_{a_{i0}}^{-1} $, where the overline indicates that all exponents in the formula for $\la(b)$ have been multiplied by $-1$. That is, $\la(b_{i1}) = \i_{g_{a_{i0}}}({\la(b)})$. Together with $g_{b_{i1}} = \i_{g_{a_{i0}}}(g_b)$, this observation motivates the following.

\begin{definition} \label{peri}
The \emph{peripheral structure} of $AC(D)$ is a set $\P(D)$ of ordered pairs of elements of $AC(D)$, defined recursively as follows.
\begin{enumerate}
    \item If $b \in AC(D)$ then $(g_b, \la(b)) \in \P(D)$. 
    \item If $(m,\la) \in \P(D)$ and $a \in A(D)$ then $(\i_{g_a}(m),\i_{g_a}(\la)) \in \P(D)$.
\end{enumerate}
\end{definition}

The elements of $\P(D)$ are the \emph{meridian-longitude pairs} in $AC(D)$. We use the notation $\P(D) = \P_1(D) \cup \dots \cup \P_\mu(D)$, where each $\P_i(D)$ includes all of the meridian-longitude pairs $(m,\la)$ with $\k_D(m)=i$. We should emphasize that the elements of a meridian-longitude pair are independent; two distinct pairs can share their meridional or longitudinal elements. See Section \ref{fourcomps} for examples.

Rather mechanical arguments using Reidemeister moves can be used to prove the following. We leave the details to the reader.

\begin{theorem} \label{invariance}
If $D$ and $D'$ are equivalent virtual link diagrams then there is an isomorphism $f:AC(D) \to AC(D')$ that maps $\P_i(D)$ to $\P_i(D')$ for each $i \in \{1, \dots, \mu\}$.
\end{theorem}

It follows that also $f(\M_i(D)) = \M_i(D')$, for each index $i$. We will see in Section \ref{im} that this implies that the subgroups $\mathcal  {I} (D) \subset \Aut(AC(D))$ and $\mathcal{I}(D') \subset \Aut(AC(D'))$ must also correspond to each other under $f$. In light of Theorem \ref{invariance}, we often replace the notation $AC(D), \M(D), \P(D), \allowbreak \mathcal{I}(D)$ with $AC(L), \M(L), \P(L)$, and $ \mathcal{I}(L)$. 

A very important algebraic invariant of an oriented link $L$ is the group $G(L)$, conventionally described by the Wirtinger presentation associated with a diagram $D$. The Wirtinger presentation has a generator $x_a$ for each $a \in A(D)$ and a relator $x_{a(c)}x_{b_1(c)}x_{a(c)}^{-1}x_{b_2(c)}^{-1}$ for each $c \in C(D)$, with $b_1(c)$ now specifically representing the underpassing arc on the right-hand side of $a(c)$ at $c$. The peripheral structure of $AC(L)$ is modeled on the classical peripheral structure of $G(L)$, with the elements of $\mathcal{I}$ playing a role analogous to the classical role of the inner automorphisms of $G(L)$.

For a classical link $L$, Boileau and Zimmerman \cite{BZ} discussed the quotient of $G(L)$ by the normal subgroup generated by the squares of meridians; this is the $\pi$-orbifold group $O(L)$. The definition extends directly to virtual links. The difference between the relators of Definition \ref{arccore} and the Wirtinger relators disappears when the squares of meridians are modded out, so there is a natural isomorphism between $O(L)$ and the quotient of $AC(L)$ modulo its normal subgroup $N$ generated by squares of meridians. The peripheral structures of $AC(L)$ and $G(L)$ both define peripheral structures on $O(L)$. These two peripheral structures on $O(L)$ have the same meridians; they differ in the choice of preferred longitudes, but it is a simple matter to obtain either family of preferred longitudes from the other. 

It turns out that when $AC(L)$ and $O(L)$ are considered as groups with specified peripheral structures, they provide equivalent link invariants. That is:

\begin{theorem} \label{piorb}
If $L_1$ and $L_2$ are virtual links, then there is an isomorphism $AC(L_1) \cong AC(L_2)$ that preserves meridian-longitude pairs if and only if there is an isomorphism $O(L_1) \cong O(L_2)$ that preserves meridian-longitude pairs.
\end{theorem}

In Section \ref{im} we discuss the involutory meridional automorphisms, and in Section \ref{piorbsec} we discuss the relationship between $AC(L)$ and $O(L)$.  In Section \ref{ab} we provide some more information about the algebraic properties of the peripheral structure of $AC(D)$, including the connection with the peripheral structure of the reduced Alexander module \cite{periadd, T1}. The longest part of the paper is Section \ref{Examples}, where we provide examples to illustrate some properties of $AC(D)$ and its peripheral structure. In Section \ref{trot} we observe that the longitudes of the core group can be used instead of the longitudes of the knot group in Trotter's proof of the noninvertibility of certain pretzel knots \cite{trotter}. In Section \ref{fourcomps} we discuss two links $L,L'$ such that $AC(L) \cong AC(L')$ but no such isomorphism matches the involutory meridional automorphisms of $L$ to those of $L'$. In Section \ref{homeo} we mention that it is possible to distinguish links with homeomorphic complements by their core groups, even without considering the peripheral structures. In Section \ref{Borr} we show that in addition to distinguishing the oriented versions of Trotter's pretzel knots, the longitudes in $AC(L)$ can be used to distinguish the oriented versions of the Borromean rings. 

The authors are grateful to Cornelius Pillen for helpful calculations. 

\section{The involutory meridional automorphisms}\label{im}

In this section we discuss some algebraic properties of these automorphisms of $AC(D)$. Most of the section is devoted to showing that $\mathcal{I}(D)$ is determined by the set $\M(D)$; it is not important which meridians are represented by arcs of a diagram $D$. First we justify Definition \ref{immap}.

\begin{lemma} \label{imlemma1}
If $D$ is a virtual link diagram and $a \in A(D)$ then $AC(D)$ has an automorphism $i_{g_a}$ with $i_{g_a}(g_b) = g_a g_b^{-1} g_a \thickspace \allowbreak \forall b \in A(D)$.
\end{lemma}
\begin{proof}
Let $F$ be the free group on the set $A(D)$, and let $a$ be a fixed element of $A(D)$. Then there is certainly an endomorphism $I_a:F \to F$ given by $I_a(b) = ab^{-1}a \thickspace \allowbreak \forall b \in A(D)$. Notice that for any $b \in A(D)$, $b = a \cdot a^{-1}ba^{-1} \cdot a=I_a(a) \cdot I_a(b)^{-1} \cdot I_a(a) = I_a(a b^{-1}a) =I_a(I_a(b))$. $F$ is generated by the elements $b \in A(D)$, so $I_a \circ I_a$ is the identity map of $F$. It follows that $I_a$ is an automorphism of $F$.

For each classical crossing $c$ of $D$, let $r'_c=a(c) b_1(c)^{-1} a(c) b_2(c)^{-1} \in F$. Then the image of $r'_c$ in the generators of $AC(D)$ is the relator $r_c$ corresponding to $c$ in Definition \ref{arccore}. Notice that
\begin{equation} \label{eq1}
\begin{aligned}
    I_a(r'_c) &= I_a(a(c)) I_a(b_1(c))^{-1} I_a(a(c)) I_a(b_2(c))^{-1} \\
    &= a a(c)^{-1} a (a b_1(c)^{-1} a)^{-1}   a a(c)^{-1} a (a b_2(c)^{-1} a)^{-1} \\
    & = a a(c)^{-1} a a ^{-1} b_1(c) a^{-1} a a(c)^{-1} a a^{-1} b_2(c) a^{-1}\\
    &= a a(c)^{-1} b_1(c)  a(c)^{-1}  b_2(c) a^{-1}\\
    &=a a(c)^{-1} b_1(c)  a(c)^{-1} \cdot (r'_c)^{-1}  \cdot(a a(c)^{-1} b_1(c)  a(c)^{-1})^{-1} \\
    &= I_a(a(c)b_1(c)^{-1} a(c) a^{-1} )\cdot (r'_c)^{-1}  \cdot I_a(a(c)b_1(c)^{-1} a(c) a^{-1} )^{-1}\text{,} 
\end{aligned}
\end{equation}
and hence 
\begin{equation} \label{eq2}
I_a \Big( (a(c)b_1(c)^{-1} a(c) a^{-1} )^{-1} \cdot r'_c \cdot a(c)b_1(c)^{-1} a(c) a^{-1} \Big) = (r'_c)^{-1} .
\end{equation}
Let $K \triangleleft F$ be the normal subgroup generated by the $r'_c$ elements. The fact that $I_a$ is an automorphism of $F$ implies that $I_a(K)$ is the normal subgroup of $F$ generated by the $I_a(r'_c)$ elements. Then (\ref{eq1}) tells us that $I_a(r'_c)$ is always a conjugate of an element of $K$, so $I_a(r'_c) \in K \thickspace \allowbreak \forall c \in C(D)$; hence $I_a(K) \subseteq K$. Also, (\ref{eq2}) tells us that $(r'_c)^{-1} \in I_a(K) \thickspace \allowbreak \forall c \in C(D)$; hence $K \subseteq I_a(K)$. Therefore $I_a(K)=K$, so $I_a$ induces an automorphism of $F/K = AC(D)$.
\end{proof}

\begin{lemma}\label{imlemma2}
If $a \in A(D)$ then $\i_{g_a}(m)=g_a m^{-1} g_a \thickspace \allowbreak \forall m \in \M(D)$. 
\end{lemma}
\begin{proof}
According to Lemma \ref{imlemma1}, $\i_{g_a}(g_b)=g_a g_b^{-1} g_a \thickspace \allowbreak \forall b \in A(D)$. 

Now, suppose $m \in \M(D)$ and the lemma holds for $m$, i.e., $\i_{g_\a}(m) = g_\a m^{-1} g_\a \thickspace \allowbreak \forall \a \in A(D)$. Then for any $\a \in A(D)$, 
\begin{align*}
\i_{g_a}(\i_{g_\a}(m)) & = \i_{g_a}(g_\a m^{-1} g_\a) = \i_{g_a}(g_\a) \i_{g_a}(m)^{-1} \i_{g_a}(g_\a) \\
& = g_a g_\a^{-1} g_a(g_a m^{-1} g_a)^{-1} g_a g_\a^{-1} g_a \\
& = g_a g_\a^{-1} g_a g_a^{-1} m g_a^{-1} g_a g_\a^{-1} g_a \\
& = g_a g_\a^{-1} m g_\a^{-1} g_a = g_a \i_{g_\a}(m)^{-1} g_a.
\end{align*}
That is, the lemma holds for $\i_{g_\a}(m)$.

The lemma follows, using induction on the number of applications of $\i_{g_\a}$ maps needed to obtain $m$ from a $g_b$ element.\end{proof}

\begin{lemma}\label{imlemma3}
Let $m_1,m_2 \in \M(D)$, and suppose there are automorphisms of $AC(D)$ given by $\i_{m_1}(n) = m_1 n^{-1} m_1$ and $\i_{m_2}(n) = m_2 n^{-1} m_2 \thickspace \allowbreak \forall n \in \M(D)$. Then there is also an automorphism $\i_{\i_{m_1}(m_2)}$ given by $\i_{\i_{m_1}(m_2)}(n) = \i_{m_1}(m_2)n^{-1} \i_{m_1}(m_2) \thickspace \allowbreak \forall n \in \M(D)$.
\end{lemma}
\begin{proof}
The composition $ \i_{m_1} \circ  \i_{m_2} \circ  \i_{m_1}$ is certainly an automorphism of $AC(D)$. For any $n \in \M(D)$,
\begin{align*}
\i_{m_1} \circ  \i_{m_2} \circ  \i_{m_1}(n) & = \i_{m_1}  \i_{m_2} (m_1 n ^{-1} m_1)\\
& = \i_{m_1} (m_2 m_1^{-1}m_2 (m_2 n^{-1} m_2)^{-1}m_2 m_1^{-1}m_2)\\
&=\i_{m_1} (m_2 m_1^{-1} n m_1^{-1}m_2)\\
&= m_1 m_2 ^{-1} m_1 (m_1 m_1^{-1} m_1)^{-1} m_1 n^{-1} m_1 (m_1 m_1^{-1} m_1)^{-1} m_1 m_2 ^{-1} m_1\\
& =m_1 m_2 ^{-1} m_1 n^{-1}m_1 m_2 ^{-1} m_1= \i_{m_1}(m_2)n^{-1} \i_{m_1}(m_2) \text{,}
\end{align*}
so $\i_{m_1} \circ  \i_{m_2} \circ  \i_{m_1}=\i_{\i_{m_1}(m_2)}$.\end{proof}

\begin{cor} \label{imcor}
For every $m \in \M(D)$ there is an automorphism $\i_m$ of $AC(D)$ with $\i_m(n) = m n^{-1} m \thickspace \allowbreak \forall n \in \M(D)$.\end{cor}
\begin{proof}
If $m=g_a$ then Lemma \ref{imlemma2} tells us that $m$ satisfies the corollary. The general case of the corollary then follows from Lemma \ref{imlemma3}, using induction on the number of applications of $\i_{g_\a}$ maps needed to obtain $m$ from a $g_a$ element. \end{proof}

Corollary \ref{imcor} tells us that we can define $\mathcal{I}(D)$ to be the subgroup of $\Aut(AC(D))$ generated by the $\i_m$ maps, $m \in \M(D)$, and each $\i_m$ map can be defined by the formula $\i_m(n) = mn^{-1}m \thickspace \allowbreak \forall n \in \mathcal{M}(D)$. Notice that we can define $\mathcal{I}(D)$ without knowing which meridians correspond to arcs in a diagram.

\begin{prop} \label{hom}
There is a homomorphism $\i:AC(D) \to \mathcal{I}(D)$ given by $\i(m)= \i_m \thickspace \allowbreak \forall m \in \mathcal{M}(D)$.
\end{prop}
\begin{proof}
In the proof of Lemma \ref{imlemma3} we showed that 
\begin{equation}\label{homeq}\i_{m_1} \circ  \i_{m_2} \circ  \i_{m_1}=\i_{\i_{m_1}(m_2)} \thickspace \allowbreak \forall m_1,m_2 \in \mathcal{M}(D).
\end{equation}
In particular, if $c \in C(D)$ has overpassing arc $a$ and underpassing arcs $b_1,b_2$ then $\i_{g_a} \circ  \i_ {g_{b_1}}  \circ  \i_{g_a} = \i_{g_{b_2}}$. The map $\i_{g_{b_1}} $ is an involution, so it equals its inverse; hence  $\i_{g_a} \circ  \i_ {g_{b_1}} ^{-1} \circ  \i_{g_a} = \i_{g_{b_2}}$. That is, the relations of Definition \ref{arccore} are satisfied by the automorphisms $\i_{g_a}$, $a \in A(D)$. It follows that there is a homomorphism $\i:AC(D) \to \mathcal{I}(D)$ with $\i(g_a) = \i_{g_a} \thickspace \allowbreak \forall a \in A(D)$.

Now, notice that if $m_1,m_2 \in \mathcal{M}(D)$ have $\i(m_1)=\i_{m_1}$ and $\i(m_2)=\i_{m_2}$ then (\ref{homeq}) implies that 
\[
\i(\i_{m_1}(m_2)) = \i(m_1 m_2^{-1} m_1)=\i(m_1)\i(m_2)^{-1} \i(m_1)= \i(m_1) \circ \i(m_2) \circ \i(m_1)= \i_{\i_{m_1}(m_2)} .
\]
The fact that $\i(m)= \i_m \thickspace \allowbreak \forall m \in \mathcal{M}(D)$ follows, using induction on the number of applications of $\i_{g_a}$ maps needed to obtain $m$ from a $g_b$ element of $AC(D)$.
\end{proof}

\begin{figure} [bht]
\centering
\begin{tikzpicture} 

\draw [thick, domain=-30:295] plot ({-6+(0.8)*cos(\x)}, {(0.8)*sin(\x)});
\draw [thick, domain=150:475] plot ({-4.95+(0.8)*cos(\x)}, {(0.8)*sin(\x)});
\node at (-6,1.1) {$a$};
\node at (-4.95,1.1) {$b$};

\end{tikzpicture}
\caption{The Hopf link, $H$.}
\label{hfig}
\end{figure}
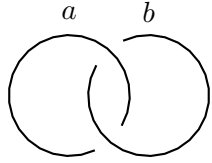

%TROTTER'S THEOREM
If $N$ denotes the normal subgroup of $AC(D)$ generated by the squares of all meridians, then the fact that $\i_m$ is always an involution implies $N \subseteq \ker \i$. As we will see in Proposition \ref{orbifold}, many links actually have $N=\ker \i$. This equality is not true for all links, however. For a simple example, consider the Hopf link $H$ pictured in Fig.\ \ref{hfig}. It has $AC(H) \cong \Z * \Z_2$, with the free product's factors generated by $g_a$ and $g_bg_a^{-1}$, respectively. Therefore the quotient of $AC(H)$ modulo the normal subgroup generated by the squared meridians is isomorphic to $\Z_2 \oplus \Z_2$. However, Definition \ref {arccore} implies that both $\i_{g_a}$ and $\i_{g_b}$ are the identity map of $AC(H)$, so $\mathcal{I}(H)$ is the trivial group.

Before proceeding, we remark that although the role of $\mathcal{I}(D)$ in the peripheral theory of core groups is analogous to the role of the inner automorphism group in the peripheral theory of link groups, there are several significant differences between inner automorphisms and the elements of $\mathcal{I}(D)$. One difference is that the inner automorphisms are determined by the structure of the link group, but $\mathcal{I}(D)$ is not determined by the group structure of $AC(D)$: determining $\mathcal{I}(D)$ requires knowing which subset of $AC(D)$ is $\M(D)$. See Section \ref{fourcomps} for examples. Another difference is that the inner automorphism given by a group element $g$ follows the same formula $x \mapsto gxg^{-1}$ for every group element $x$, but an $i_m$ map is not given by the same formula for all group elements. For instance, $\i_{g_a}(g_b^{-1})$ must equal the inverse of $i_{g_a}(g_b) = g_ag_b^{-1} g_a$; this is not the same as $g_a (g_b^{-1})^{-1} g_a$.

\section{The $\pi$-orbifold group} \label{piorbsec}

Let $V$ be the 2-fold cover of $\S^3$ branched over a classical link $L$, and let $\tilde V$ be its universal cover. In \cite{BZ} Boileau and Zimmermann define the $\pi$-\emph{orbifold group} of $L$, denoted by $O(L)$, to be the group of diffeomorphisms of $\tilde V$ generated by involutions such that the quotient space $\tilde V/O(L)$ is an orbifold with underlying space $\S^3$ and branch locus $L$. They show that $O(L)$ is isomorphic to the link group $\pi_1(\S^3 \setminus L)$ modulo the normal subgroup generated by squares of the meridians of $L$. This result allows us to extend the definition of $O(L)$ to virtuals: if $L$ is a virtual link with group $G(L)$ described by the Wirtinger presentation, then $O(L)$ is the quotient of $G(L)$ modulo the normal subgroup of $G(L)$ generated by squares of meridians.

As mentioned in the introduction, a Wirtinger relator at a crossing of a link diagram $D$ is the same as the core group relator when the additional relators $g_a^2, a \in A(D)$ are added to the presentation of $AC(D)$ in Definition \ref{arccore}. Hence $O(L)$ is also the quotient of the arc core group of $L$ modulo the normal subgroup $N$ generated by the squares of meridians.

In \cite{BZ} a prime, unsplittable link in $\S^3$ with an infinite $\pi$-orbifold group is said to be \emph{sufficiently complicated} (abbv.\ s.c.). 
Several deep results about s.c.\ links are proven in \cite{BZ}. Theorem 1 of \cite{BZ} states that if $L, L'$ are s.c.\ links in $\S^3$, then $(\S^3, L)$ and $(\S^3, L')$ are diffeomorphic if and only if their $\pi$-orbifold groups are isomorphic. Theorem 3.1 of \cite{BZ} extends this result to cover all prime, unsplittable links that are not 2-bridge links.

Let ${\rm Sym}(\S^3, L)$ denote the symmetry group of a link $L$ in $\S^3$, the group of pairwise isotopy classes of diffeomorphisms of $(\S^3, L)$. Diffeomorphisms and isotopies can be lifted to $\tilde V$, and hence there is a homomorphism $\gamma$ from  ${\rm Sym}(\S^3, L)$ to the outer automorphism group  ${\rm Out} (O(L))$, as lifted diffeomorphisms induce automorphisms of $O(L)$. Theorem 2 of \cite{BZ} states that if $L$ is s.c., then $\gamma$ is an isomorphism. 

Another interesting result of Boileau and Zimmerman \cite {BZ} is their Corollary 1.4, which tells us that if $L$ is s.c.\ and not Seifert-fibered then the center of $O(L)$ is trivial. Here is a consequence of this result.

\begin{prop} \label{orbifold} Assume that a classical link $L$ is s.c.\ and not Seifert-fibered. Then the involutory meridional automorphism group $\mathcal{I}(L)$ is isomorphic to the $\pi$-orbifold group $O(L)$. 
\end{prop}

\begin{proof} By Proposition \ref{hom} there is a homomorphism $\iota: AC(L) \to \mathcal{I}(L)$  sending $m \in \mathcal{M}(L)$ to $\iota_m$. Since $\iota$ is trivial on squares of meridians, it induces a homomorphism 
$\bar \iota : AC(L)/N \cong O(L) \to \mathcal{I}(L)$. 

Every automorphism $\iota_m \in \mathcal{I}(L)$ induces an inner automorphism $i_m$ of $O(L)$, since the quotient classes of a meridian and its inverse are the same in $O(L)$.  Moreover, the assignment $\iota_m \mapsto i_m$ defines a homomorphism  $i$ from $\mathcal{I}(L)$ to $\textup{Inn}(O(L))$, the group of inner automorphisms of $O(L)$. 

Consider the composition $O(L) \xrightarrow{\bar \iota} \mathcal{I}(L)  \xrightarrow{i} \textup{Inn}(O(L))$.
By Corollary 1.4 of \cite{BZ}, the center of $O(L)$ is trivial, and so the composition $i \circ \bar \iota$ is an isomorphism. Hence $\mathcal{I}(L) \cong O(L).$\end{proof} 

Examples of knots satisfying the hypotheses of Proposition \ref{orbifold} are found in Section \ref{trot}. The discussion of the Hopf link $H$ at the end of Section \ref{im} indicates that $H$ does not satisfy either the hypotheses or the conclusion of Proposition \ref{orbifold}.

\begin{figure} [bht]
\centering
\begin{tikzpicture} 

\draw [thick] [->] (4,-1) -- (2.4,.6);
\draw [thick] (2,1) -- (2.4,.6);
\draw [thick] (2,-1) -- (2.8,-.2);
\draw [thick] [->] (4,1) -- (3.6,.6);(4,-1)
\draw [thick] (3.2,.2) -- (3.6,.6);
\node at (1.5,-.65) {$b_{i(j+1)}$};
\node at (4.3,.65) {$b_{ij}$};
\node at (4.3,-.65) {$a_{ij}$};

\draw [thick] [<-] (-2.4,-.6) -- (-4,1);
\draw [thick] (-2.4,-.6) -- (-2,-1);
\draw [thick] (-4,-1) -- (-3.2,-.2);
\draw [thick] [->] (-2,1) -- (-2.4,.6);
\draw [thick] (-2.8,.2) -- (-2.4,.6);
\node at (-4.5,-.65) {$b_{i(j+1)}$};
\node at (-1.7,.65) {$b_{ij}$};
\node at (-1.7,-.65) {$a_{ij}$};

\node at (-3,-1.8) {$w_{ij}=-1$};
\node at (3.1,-1.8) {$w_{ij}=1$};

\end{tikzpicture}
\caption{The crossing $c_{ij}$ has writhe $w_{ij}$.}
\label{crossfig}
\end{figure}
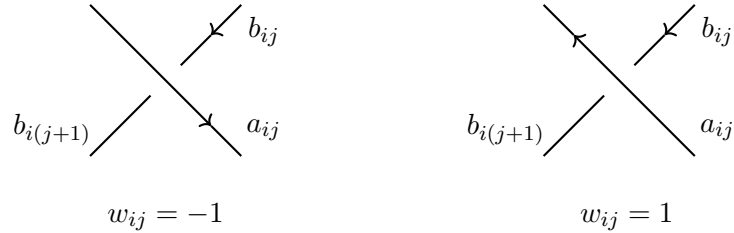

We now turn our attention to Theorem \ref{piorb}. To prepare, we establish notation for the peripheral structures of the group $G(D)$ and the $\pi$-orbifold group $O(D)$ of a virtual link diagram $D$. As mentioned in the introduction, $G(D)$ is the group given by the presentation 
\[\langle \{x_a \mid a \in A(D)\};\{x_{a(c)}x_{b_1(c)}x_{a(c)}^{-1}x_{b_2(c)}^{-1}\mid c \in C(D)\} \rangle\text{,}\]
where $b_1(c)$ is the underpassing arc on the right-hand side of the overpassing arc $a(c)$ of a crossing $c$, and $b_2(c)$ is the underpassing arc on the left of $a(c)$. The meridians of $G(D)$ are the elements $x_a$, and their conjugates. Let $b$ be an arc of $K_i$ in $D$, and index the arcs of $K_i$ in $D$ as in Definition \ref{peri}. That is, the arcs of $K_i$ are indexed as $b=b_{i0},b_{i1},\dots, b_{i(k_i-1)},b_{ik_i}=b$, following the orientation of $K_i$. Also, for $0 \leq j \leq k_i$ the crossing of $D$ separating $b_{ij}$ from $b_{i(j+1)}$ is $c_{ij}$, and $a_{ij}$ is the overpassing arc at $c_{ij}$. Then the preferred longitude at $b$ is conventionally defined to be 
\begin{equation*} 
x_{a_{i0}}^{w_{i0}}x_{a_{i1}}^{w_{i1}} \dots x_{a_{i(k_i-1)}}^{w_{i(k_i-1)}} x_b^p\text{,}
\end{equation*}
where $w_{ij}$ is the writhe of $c_{ij}$ (see Figure \ref{crossfig}) and $p \in \Z$ is chosen so that the total exponent with respect to arcs of $K_i$ is $0$.

The $\pi$-orbifold group $O(D)$ can be defined to be the quotient of $G(D)$ modulo the normal subgroup generated by the squares of the $x_a$ elements, and the peripheral structure of $O(D)$ is the image of the peripheral structure of $G(D)$ under the canonical surjection $\s:G(D) \to O(D)$. Notice that the only part of the peripheral structure of $O(D)$ that is sensitive to the orientations of link components is the choice of preferred longitudes. As mentioned above, it is clear that if $N$ is the normal subgroup of $AC(D)$ generated by the squares of the $g_a$ elements then there is a natural isomorphism $O(D) \cong AC(D)/N$, with $\s(x_a)$ corresponding to $g_aN$ for every $a \in A(D)$. We use this isomorphism to identify the two groups. The $\i_m$ maps induce inner automorphisms of $AC(D)/N$, because the difference between exponents $+1$ and $-1$ disappears in the quotient, so in general, meridians of $AC(D)$ are mapped to meridians of $O(D)$.

If $b \in A(D)$ then the image of $\la(b)$ in $AC(D)/N$ is not necessarily the same as the image under $\s$ of the conventional preferred longitude of $G(D)$. Both images are of the form 
\[
g_{a_{i0}}g_{a_{i1}} \dots g_{a_{i(k_i-1)}}g_b^q N \text{.}
\]
For the image of $\la (b)$, $q$ makes the total exponent even, while for the image of the preferred longitude of $G(D)$, $q$ makes the total exponent with respect to arcs of $K_i$ even. Although they are not the same, either longitude is easily obtained from the other. It follows that the peripheral structure of $O(D)$ is determined directly by the peripheral structure of $AC(D)$.

To complete the proof of Theorem \ref{piorb}, we show that the peripheral structure of $O(D)$ determines the peripheral structure of $AC(D)$. An essential ingredient of the proof is Theorem \ref{lplus} below, which was first mentioned for classical links by Sakuma \cite{S}; later accounts have been given by Pryzytcki \cite{P} and Silver, Traldi and Williams \cite{STW1}. To set up the statement, let $D$ be a virtual link diagram, and let $D^+$ be the diagram obtained from $D$ by inserting a new unknotted component, which does not appear in any crossing. The one arc of this new component is denoted $a^+$, and the canonical surjection $G(D^+) \to O(D^+)$ is denoted $\s^+$. 

\begin{theorem} \label{lplus}
There is a monomorphism $f:AC(D) \to O(D^+)$ given by $f(g_a) = \s^+(x_a x_{a^+}) \thickspace \allowbreak \forall a \in A(D)$. The image of $f$ is the subgroup of $O(D^+)$ consisting of products of even numbers of meridians.
\end{theorem}

We refer to \cite{P, S, STW1} for full proofs of Theorem \ref{lplus}, but we take a moment to discuss the description of the image of $f$. Suppose $k \in \N$ is an even integer, $a_1, \dots, a_k \in A(D^+)$ and $g=\s^+(x_{a_1} \dots x_{a_k}) \in O(D^+)$. We may presume that no two consecutive $a_i$ are equal, because $\s^+(x_a^2)=1$ for every arc $a \in A(D^+)$. If $i \geq 1$ and $a_{2i-1},a_{2i}$ are both arcs of $D$, then $\s^+(x_{a_{2i-1}}x_{a_{2i}})= \s^+(x_{a_{2i-1}} x^2_{a^+}x_{a_{2i}})=f(g_{a_{2i-1}}g^{-1}_{a_{2i}})$.  If $i \geq 1$ and $a_{2i-1}=a^+$ then $a_{2i}$ is an arc of $D$, and $\s^+(x_{a_{2i-1}}x_{a_{2i}}) = \s^+(x_{a^+}x_{a_{2i}})=f(g^{-1}_{a_{2i}})$. If $i \geq 1$ and $a_{2i}=a^+$ then $a_{2i-1}$ is an arc of $D$, and $\s^+(x_{a_{2i-1}}x_{a_{2i}}) = \s^+(x_{a_{2i-1}}x_{a^+})=f(g_{a_{2i-1}})$. Therefore if  $b_1, \dots, b_j$ is the sublist of $a_1, \dots, a_k$ obtained by removing all occurrences of $a^+$, then $g = f(g_{b_1}^{\e_1} g_{b_2}^{\e_2} \dots g_{b_j}^{\e_j})$ for some choice of $\e_1, \dots, \e_j \in \{-1,1\}$. To be precise, $\e_i = (-1)^{i'-1}$, where $b_i$ is $a_{i'}$.

We will need some more notation regarding $D^+$. There are natural inclusion maps $AC(D) \subset AC(D^+)$ and $O(D) \subset O(D^+)$, which we leave implicit; i.e., we do not adopt any notation for them. Then $AC(D^+)$ is the internal free product $AC(D) * \Z$, with the second free factor generated by $g_{a^+}$, and $O(D^+)$ is the internal free product $O(D)*\Z_2$, with the second free factor generated by $\s^+(x_{a^+})$. There are canonical projections onto the quotients, $\pi_A:AC(D^+) \to AC(D)$ and $\pi_O:O(D^+) \to O(D)$; $\ker \pi_A$ is the normal closure of $g_{a^+}$ and $\ker \pi_O$ is the normal closure of $\s^+(x_{a^+})$. Let $N^+$ be the normal subgroup of $AC(D^+)$ generated by the squares of meridians. Then as discussed above for $AC(D)$ and $O(D)$, there is an isomorphism $O(D^+) \cong AC(D^+)/N^+$, under which $\s^+(x_a)$ corresponds to $g_a N^+$ for each $a \in A(D^+)$. We use this isomorphism to identify $AC(D^+)/N^+$ and $O(D^+)$. Notice that the composition $AC(D) \subset AC(D^+) \to AC(D^+)/N^+ \cong O(D^+)$ is a homomorphism $AC(D) \to O(D^+)$ denoted by the deceptively simple formula $x \mapsto xN^+$; unlike $f$, this homomorphism is not injective in general.

\begin{lemma} \label{lpluslem}
The monomorphism $f$ of Theorem \ref{lplus} has $f(m) = (m N^+) \s^+(x_{a^+}) \thickspace \allowbreak \forall m \in {\cal M}(D)$. Moreover, if $m \in {\cal M}(D^+)$ then $f^{-1}((mN^+) \s^+(x_{a^+})) \in {\cal M}(D)$ if and only if $m N^+=\pi_O(m N^+)$.
\end{lemma}
\begin{proof}
The first sentence of the lemma holds if $m=g_a$ for some $a \in A(D)$, because $f(g_a) = \s^+(x_a) \s^+(x_{a^+})=(g_a N^+) \s^+(x_{a^+})$.

Suppose the first sentence of the lemma holds for $m,n\in {\cal M}(D)$. Then 
\begin{align*}
f(\i_n(m))&=f(nm^{-1}n)=f(n)f(m)^{-1}f(n)\\
&=(n N^+) \s^+(x_{a^+}) ((m N^+) \s^+(x_{a^+}))^{-1} (n N^+) \s^+(x_{a^+})\\
&= (n N^+) \s^+(x_{a^+}) (\s^+(x_{a^+}))^{-1} (m^{-1}N^+) (n N^+) \s^+(x_{a^+})\\
&= (n N^+) (m^{-1}N^+) (n N^+) \s^+(x_{a^+})\\
&= (n m^{-1} n N^+) \s^+(x_{a^+})= (\i_n (m) N^+) \s^+(x_{a^+}) \text{,}
\end{align*}
so the first sentence of the lemma holds for $\i_n(m)$ too. Using induction on the number of applications of $\i_{g_b}$ maps needed to obtain $m$ from a $g_a$ element, it follows that the first sentence of the lemma holds for every $m \in {\cal M}(D)$. 

The second sentence of the lemma is not so easy to prove.

If $m$ is a meridian of $AC(D^+)$ associated with the new component then when $m$ is written as a product of $g_a$ elements, the exponent sum with respect to arcs of $D$ is even; hence the exponent sum of $\pi_A(m)$ is even. The exponent sum of a meridian is odd, so $\pi_A(m) N^+ = \pi_O(m N^+)$ cannot equal $m N^+$. According to the discussion after Theorem \ref{lplus}, $f^{-1}((m N^+) \s^+(x_{a^+}))$ will also have an even exponent sum, so  $f^{-1}((m N^+)\s^+(x_{a^+}))$ cannot be a meridian in $AC(D)$.

Now, suppose $m$ is a meridian of $AC(D^+)$ associated with a component of $D$. Then $mN^+$ is a meridian of $O(D^+)$ associated with the same component of $D$, so it is a conjugate of $\s^+(x_a)$, for some $a \in A(D)$. Say $m N^+=g\s^+(x_a)g^{-1}$, where $g \in O(D^+)$. Replacing $g$ with $g \s^+(x_a)$ if necessary, we may presume that $g=\s^+(x_{a_1} \dots x_{a_k})$ for some even integer $k \in \N$ and some $a_1, \dots, a_k \in A(D^+)$.  As discussed after Theorem \ref{lplus}, it follows that $g = f(g_{b_1}^{\e_1} g_{b_2}^{\e_2} \dots g_{b_j}^{\e_j})$ for some choice of $b_1, \dots, b_j \in A(D)$ and $\e_1, \dots, \e_j \in \{-1,1\}.$

The formula $f(g_b) = \s^+(x_b x_{a^+})$ implies that for every $b \in A(D)$, we have these two equalities.
\begin{align*}
&\s^+(x_{a^+})f(g_b)= \s^+(x_{a^+}x_b) \s^+(x_{a^+})=f(g_b^{-1})\s^+(x_{a^+}) \\
&\s^+(x_{a^+})f(g_b^{-1})= \s^+(x_b)=f(g_b)\s^+(x_{a^+}) \end{align*}
Using these equalities, we have 
\begin{align*}
m N^+ &=g\s^+(x_a)g^{-1} \\
&= f(g_{b_1}^{\e_1} g_{b_2}^{\e_2} \dots g_{b_j}^{\e_j})\s^+(x_a)\s^+(x_{a^+}) \cdot \s^+(x_{a^+}) f(g_{b_j}^{-\e_j} \dots g_{b_2}^{-\e_2} g_{b_1}^{-\e_1})\\
&=f(g_{b_1}^{\e_1} g_{b_2}^{\e_2} \dots g_{b_j}^{\e_j}) f(g_a) f(g_{b_j}^{\e_j} \dots g_{b_2}^{\e_2} g_{b_1}^{\e_1})\s^+(x_{a^+}) 
\end{align*}
and hence
\begin{equation}\label{eq3}
f^{-1}((m N^+) \s^+(x_{a^+}) ) = g_{b_1}^{\e_1} g_{b_2}^{\e_2} \dots g_{b_j}^{\e_j} g_a g_{b_j}^{\e_j} \dots g_{b_2}^{\e_2} g_{b_1}^{\e_1}.
\end{equation}

Suppose that $m N^+ = \pi_O(m N^+)$. Then we can replace $g$ with $\pi_O(g)$ in the equation $m N^+ = g \s^+(x_a) g ^{-1}$. Then $a_1, \dots, a_k$ are all arcs of $D$. Referring again to the discussion after Theorem \ref{lplus}, it follows $j=k$, $a_i = b_i$ for every index $i$, and $\e_i = (-1)^{i-1}$ for each index $i$. Then (\ref{eq3}) implies
\begin{equation*}
f^{-1}((m N^+) \s^+(x_{a^+}) )= g_{b_1}g_{b_2}^{-1} \dots g_{b_j}^{-1} g_a g_{b_j}^{-1} \dots g_{b_2}^{-1} g_{b_1} = \i_{g_{b_1}} \i_{g_{b_2}} \dots \i_{g_{b_j}}(g_a) \text{,}
\end{equation*}
which is a meridian of $AC(D)$.

Now, suppose that $f^{-1}((m N^+) \s^+(x_{a^+}) )$ is a meridian of $AC(D)$. Then (\ref{eq3}) gives an expression of $f^{-1}((m N^+) \s^+(x_{a^+}) )$ in which $g_a$ is the only meridian with an odd total exponent, so $f^{-1}((m N^+) \s^+(x_{a^+}) )$ must be a meridian associated with the same component as $g_a$. That is, there is an automorphism $\i \in {\cal I}(D)$ with $f^{-1}((m N^+) \s^+(x_{a^+}) ) = \i(g_a)$. Replacing $\i$ with $\i \circ \i_{g_a}$ if necessary, we may presume that there exist an even integer $n$ and arcs $\a_1, \dots, \a_n \in A(D)$ with $\i = \i_{g_{\a_1}} \circ \dots \circ \i_{g_{\a_n}}$. Then 
\begin{equation}\label{eq4}
\begin{aligned}
m N^+ &= ff^{-1}((m N^+) \s^+(x_{a^+}))\s^+(x_{a^+}) = f(\i(g_a)) \s^+(x_{a^+}) \\
&=f(g_{\a_1}g^{-1}_{\a_2} \dots g^{-1}_{\a_n} g_a g^{-1}_{\a_n} g_{\a_{n-1}} \dots g_{\a_1})\s^+(x_{a^+}).
\end{aligned}
\end{equation}

Notice that 
\begin{align*}    
f(g_{\a_1}g^{-1}_{\a_2}) &= \s^+(x_{\a_1}x_{a^+})\s^+(x_{a^+}x_{\a_2})\\&= \s^+(x_{\a_1})\s^+(x_{a^+})^2 \s^+(x_{\a_2}) = \s^+(x_{\a_1})\s^+(x_{\a_2}) = \s^+(x_{\a_1}
x_{\a_2})\text{,}
\end{align*}
because $\s^+(x_{a^+})^2=1$. Applying the same kind of cancellation $n$ times, (\ref{eq4}) yields
\begin{align*}
m N^+ &= f(g_{\a_1}g^{-1}_{\a_2} \dots g^{-1}_{\a_n} g_a g^{-1}_{\a_n} g_{\a_{n-1}} \dots g_{\a_1})\s^+(x_{a^+})\\
&= f(g_{\a_1}g^{-1}_{\a_2}) \dots f(g_{\a_{n-1}}g^{-1}_{\a_n}) f(g_a g^{-1}_{\a_n}) \dots  f(g_{a_3} g_{\a_2}^{-1})f(g_{\a_1})\s^+(x_{a^+})\\
&=\s^+(x_{\a_1}x_{\a_2} \dots x_{\a_n} x_a x_{\a_n}  \dots x_{\a_3}x_{\a_2}) \s^+(x_{\a_1} ) \s^+(x_{a^+})\s^+(x_{a^+})\\
&= \s^+(x_{\a_1}x_{\a_2} \dots x_{\a_n} x_a x_{\a_n}  \dots x_{\a_1} )\s^+(x_{a^+})^2\\
&= \s^+(x_{\a_1}x_{\a_2} \dots x_{\a_n} x_a x_{\a_n} x_{\a_{n-1}} \dots x_{\a_1} ).
\end{align*}
As $\a_1,\dots,\a_n$ and $a$ are all arcs of $D$, it follows that 
\[m N^+ = \pi_O(\s^+(x_{\a_1}x_{\a_2} \dots x_{\a_n} x_a x_{\a_n} x_{\a_{n-1}} \dots x_{\a_1} ))=\pi_O(m N^+ ). \qedhere 
\]
\end{proof}

Using Theorem \ref{lplus} and Lemma \ref{lpluslem}, it is not hard to show that the peripheral structure of $O(D)$ determines the peripheral structure of $AC(D)$ in two steps.

1. The peripheral structure of $O(D)$ determines both the peripheral structure of $O(D^+)$ and the element $\s^+(x_{a^+})$, as follows. The group $O(D^+)$ is the free product $O(D)*\Z_2$, with the free factor $\Z_2$ generated by $\s^+(x_{a^+})$. The peripheral structure of $O(D^+)$ includes two kinds of meridian-longitude pairs: if $g \in O(D^+)$ and $(m, \la)$ is a meridian-longitude pair of $O(D)$ then $(gmg^{-1},g \la g^{-1})$ is a meridian-longitude pair of $O(D^+)$, and if $g \in O(D^+)$ then $(g \s^+(x_{a^+})g^{-1},1)$ is a meridian-longitude pair of $O(D^+)$.

2. The peripheral structure of $O(D^+)$, together with the element $\s^+(x_{a^+})$, determines the peripheral structure of $AC(D)$, as follows. The meridians of $AC(D)$ are the elements $f^{-1}(m \s^+(x_{a^+}))$, where $m$ is a meridian of $O(D^+)$ that has $m=\pi_O(m)$. If $m$ is a meridian of $O(D^+)$ with $m=\pi_O(m)$ and $(m,\la)$ is a meridian-longitude pair of $O(D^+)$, let $q \in \{0,1\}$ be congruent (mod 2) to the exponent sum of a product of meridians of $O(D^+)$ equal to $\la$. Then $(f^{-1}(m \s^+(x_{a^+})), f^{-1}(\la m^q))$ is a meridian-longitude pair of $AC(D)$.

This completes the proof of Theorem \ref{piorb}. Notice that the proof shows that for unoriented links, the combination of $AC(L)$ and $\M(L)$ is equivalent to the combination of $O(L)$ and its set of meridians. As a consequence, we deduce from Theorem 3.1 of \cite{BZ} that up to mirror images, the combination of $AC(L)$ and $\M(L)$ is a classifying invariant for unoriented prime, unsplittable links that are not 2-bridge links. That is, two such links are equivalent up to mirror images if and only if there is a meridian-preserving isomorphism between their $AC$ groups.

\section{Some other properties of $AC(L)$} \label{ab}

In this section we briefly summarize some properties of meridians and longitudes in core groups. 

\begin{prop} \label{i1lemma}
A core group $AC(L)$ has an automorphism $\i_1$ with $\i_1(m) = m^{-1} \thickspace \allowbreak \forall m \in \mathcal{M}(L)$.
\end{prop}
\begin{proof}
First, we claim that $AC(L)$ has an automorphism $\i_1$ with $\i_1(g_b) = g_b^{-1}  \thickspace \allowbreak \forall b \in A(D)$. The proof of the claim is a simplified version of the proof of Lemma \ref{imlemma1}, obtained by replacing $I_a$ with the automorphism $I_1:F \to F$ given by $I_1(b) = b^{-1} \thickspace \allowbreak \forall b \in A(D)$. 

Second, we claim that $\i_1(m) = m^{-1} \thickspace \allowbreak \forall m \in \mathcal{M}(L)$. This second claim holds when $m=g_b$, of course, and if the claim holds for $m,n \in \mathcal{M}(L)$ then it also holds for $\i_n(m)$, because 
\[
\i_1(\i_n(m)) = \i_1(n m^{-1} n) = \i_1( n)\i_1( m)^{-1} \i_1( n) = n^{-1} (m^{-1})^{-1} n^{-1} = (n m^{-1} n) ^{-1}=\i_n(m)^{-1}.
\]
The second claim follows, using induction on the number of applications of $\i_{g_a}$ maps needed to obtain $m$ from a $g_b$ element.
\end{proof}

\begin{prop}
\label{pm}
A core group $AC(L)$ has a subgroup
\[
AC^{\pm}(L) = \{m_1 m_2^{-1} \dots m_{2k-1}m_{2k}^{-1} \mid k\in \N \text{ and } m_1, \dots , m_{2k} \in \mathcal{M}(L) \}.
\]
If $m_0$ is any meridian of $AC(L)$ then $\i_{m_0}(AC^\pm(L)) = AC^\pm(L)$. Moreover, the cyclic subgroup $(m_0)$ is infinite, and $AC(L)$ is the internal free product $ AC^ \pm (L)*(m_0)$. 
\end{prop}
\begin{proof}
The first sentence of the statement follows from the fact that $AC^\pm(L)$ is closed under inverses and products. The second sentence follows from the fact that for any $m_1, \dots, m_{2k} \in \mathcal{M}(L)$,
\[
\i_{m_0}(m_1 m_2^{-1} \dots m_{2k-1}m_{2k}^{-1})=m_0 m_1^{-1} m_2 \dots m_{2k-1}^{-1} m_{2k} m_0^{-1} \in AC^\pm(L)
\]
and
\[
m_1 m_2^{-1} \dots m_{2k-1}m_{2k}^{-1} = \i_{m_0}(m_0 m_1^{-1} m_2 \dots m_{2k-1}^{-1} m_{2k} m_0^{-1}) \in \i_{m_0}(AC^\pm(L)).
\]

For the third sentence of the statement, assume first that $m_0 = g_{b_0}$. For each $a \in A(D)$, let $h_a = g_a g_{b_0}^{-1}$; in particular, $h_{b_0}=1$. We can rewrite Definition \ref{arccore} in terms of the generating set $\{g_{b_0}\} \cup \{ h_a \mid a \in A(D) \}$ by changing each appearance of $g_a$ in a relator to $h_a$, and inserting the relator $h_{b_0}$. No relator of this rewritten presentation involves $g_{b_0}$, so the cyclic subgroup $(g_{b_0})$ is infinite and $AC(L) \cong H* (g_{b_0})$, where $H$ is the subgroup generated by $\{h_a \mid a \in A(D)\}$. It is obvious that $H \subseteq AC^\pm(L)$, and the inclusion $AC^\pm(L) \subseteq H$ follows from the fact that a product $m_1 m_2^{-1} \dots m_{2k-1}m_{2k}^{-1}$ is equal to a product of powers of $g_a$ elements, with alternating exponents $1, -1, 1, -1, \dots$. Such a product is equal to the element of $H$ obtained by replacing each appearance of a $g_a$ element in the product by $h_a = g_a g_{b_0}^{-1}$; the powers of $g_{b_0}$ all cancel out.

Observe that if $n \in \mathcal{M}(L)$ and the third sentence of the statement holds for a meridian $m_0$, then the fact that $\i_n$ is an automorphism of $AC(L)$ with $\i_n(AC^\pm(L)) = AC^\pm (L)$ guarantees that the cyclic subgroup $(\i_n(m_0))$ is  infinite, and $AC(L)$ is the internal free product $AC^\pm(L) * (\i_n(m_0))$. That is, the third sentence of the statement holds for $\i_n(m_0)$. We deduce that the third sentence of the statement holds in general by using induction on the number of applications of $\i_n$ maps needed to obtain $m_0$ from a $g_{b_0}$ element.
\end{proof}

Notice that the longitudes of $AC(L)$ are all elements of $AC^{\pm}(L)$. Also, $AC(L)$ has an analogous subgroup $AC^\mp(L)$, with analogous properties. The automorphism $\i_1$ has $\i_1(AC^\pm(L))=AC^\mp(L)$, so $AC^\pm(L)$ and $AC^\mp(L)$ are isomorphic to each other, and equivalent as subgroups of $AC(L)$.

The subgroup we denote $AC^\pm(L)$ features prominently in the classical theory of $AC(L)$: if $L$ is a classical link then Sakuma \cite{S} and Wada \cite{W} showed that $AC^\pm(L) \cong \pi_1(X_2)$, the fundamental group of the 2-fold cover of $\S^3$ branched along $L$. Of course it follows that the abelianization $AC^\pm(L)_{ab}$ is isomorphic to the homology group $H_1(X_2)$. 

Here is a familiar property of link groups: if $L=K_1 \cup \dots \cup K_ \mu $ then $G(L \setminus K_\mu)$ is isomorphic to the quotient of $G(L)$ modulo the normal subgroup generated by a meridian of $K_\mu$. In contrast, obtaining $AC(L \setminus K_\mu)$ from $AC(L)$ requires two steps: first restrict to the subgroup $S \subset AC(L)$ generated by the meridians of $K_1, \dots , K_{\mu-1}$, and then take the quotient of $S$ modulo its normal subgroup generated by all elements of the form $\i_m(s)s^{-1}$, where $m$ is a meridian of $K_\mu$ and $s \in S$. Notice in particular that determining $AC(L \setminus K_\mu)$ involves all the meridians of $AC(D)$, not just a meridian of $K_\mu$. An example is discussed in Section \ref{fourcomps}.

The \emph{reduced Alexander module} of a link $L$ is a famous invariant; we denote it $\Mr(L)$. It is a module over the ring $\Lambda = \Z[t,t^{-1}]$ of Laurent polynomials in the variable $t$, with integer coefficients. Like Definition \ref{arccore}, the definition of $\Mr(L)$ involves generators corresponding to diagram arcs, and relations corresponding to crossings. There is a module generator $\g_a$ for each $a \in A(D)$. If $c$ is a crossing with overpassing arc $a$ and underpassing arcs $b_1,b_2$, with $b_1$ on the right-hand side of $a$, then there is a corresponding relation $\g_{b_2} = (1-t)\g_a + t \g_{b_1}$. Peripheral elements in reduced Alexander modules have been studied recently \cite{periadd, T1}.

The connection between core groups and reduced Alexander modules is simple. Let $\nu:\Lambda \to \Z$ be the homomorphism with $\nu(t)=-1$, and let $M_A(L)_\nu$ denote the abelian group obtained by applying $\nu$ to the definition of  $\Mr(L)$. Then there is a natural isomorphism between $M_A(L)_\nu$ and the abelianization of $AC(L)$; by ``natural'' we mean that the image of the group generator $g_a$ in $AC(L)_{ab}$ is the same as the image of the module generator $\g_a$ in $M_A(L)_\nu$. The longitudes in $AC(L)$ and $\Mr(L)$ also have the same image in $ AC(L)_{ab} \cong M_A(L)_\nu$. Results from \cite{periadd, mvaq2,T1, mvaq4, mvaq5} give us the following properties of these longitudes.

\begin{itemize}

\item All the meridians of a single component of $L$ have the same longitude in $AC(L)_{ab}$.

\item The longitudes in $AC(L)_{ab}$ are all elements of order 1 or 2, and they generate the subgroup of $AC(L)_{ab}$  annihilated by 2.

\item For any virtual knot $K$, the longitude is 0 in $AC(K)_{ab}$.

\item For a classical link or a checkerboard colorable virtual link of $\mu$ components, the sum of the  $\mu$ longitudes is 0 in $AC(L)_{ab}$.
\end{itemize}

Examples discussed in Section \ref{Examples} show that the longitudes in $AC(L)$ do not share these properties. A single component of a link may have several distinct associated longitudes, longitudes may be of infinite order in $AC(L)$, and a knot may have nontrivial longitudes.

\section{Examples} \label{Examples}

\subsection{Trotter's non-invertible pretzel knots} \label{trot}

\begin{figure}[H]
\begin{center}
\includegraphics[height=1.8 in]{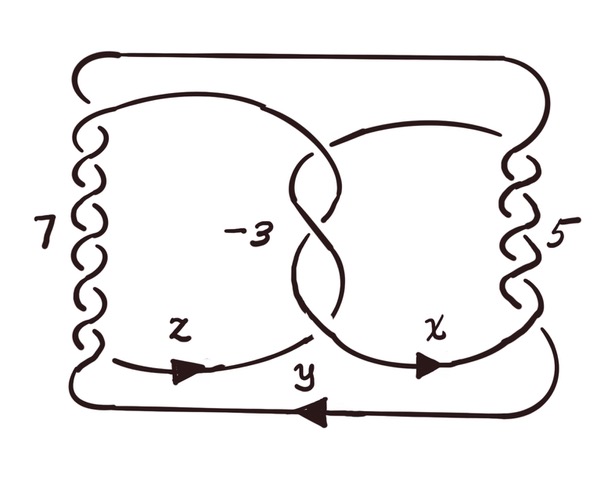}
\caption{Diagram $D$ for K(7,-3,5)}
\label{Pretzeldiagram}
\end{center}
\end{figure}

A knot is {\it invertible} if, given an orientation, it can be deformed by isotopy into itself with the orientation reversed. In \cite{trotter} Trotter gave the first examples of non-invertible knots, pretzel knots $K=K(p,q,r)$ given by odd integers $p, q, r$ such that $|p|, |q|$ and $|r|$ are distinct and greater than 1. A diagram of $K(7,-3,5)$ appears in Figure \ref{Pretzeldiagram}. It has been assigned an orientation as well as Wirtinger generators $x, y, z$ for $G=\pi_1(\S^3 \setminus K)$. (In this section we adopt the notation of \cite{trotter}, using the same letters for arcs of $D$ as for their corresponding Wirtinger generators.)

As we show, the information needed for Trotter's proof is found in the peripheral structure of the core group of $K$. 

Define $k, l, m$ by $p = 2k+1, q= 2l+1, r=2m+1$. Trotter begins his argument with a routine computation of a Wirtinger presentation for $G$. It has three generators $x, y, z$ and three relations. (As usual, any one relation is a consequence of the other two): 

$$(xy^{-1})^mx(xy^{-1})^{-m} = (yz^{-1})^{k+1}z(yz^{-1})^{-k-1}$$
$$(yz^{-1})^ky(yz^{-1})^{-k} = (zx^{-1})^{l+1}x(zx^{-1})^{-l-1}$$
$$(zx^{-1})^lz(zx^{-1})^{-l} = (xy^{-1})^{m+1}y(xy^{-1})^{-m-1}$$

A longitude $w$ commuting with $x$ is described by 

$$w= (xy^{-1})^{-m}(yz^{-1})^{k+1}(zx^{-1})^{-1}(xy^{-1})^{m+1}(yz^{-1})^{-k}(zx^{-1})^{l+1}.$$

In order for $K$ to be invertible, there must be an {\it inversion} of $G$, an automorphism of $G$ that inverts $x, y, z$ and $w$. Trotter shows that an inversion cannot exist. His {\bf first step} is to form a quotient group $G/H$, where $H$ is the normal subgroup of $G$ generated by squares of meridians. After  simplification, we have

$$G/H = \< x, y, z ; (xy)^r = (yz)^p = (zx)^q \>.$$

The subgroup $U$ generated by $(xy)^r$ is the center of the commutator subgroup of $G/H$.  
In his {\bf second step}, Trotter forms $W = (G/H)/U$. We have

$$W= \< x, y, z ; (xy)^r= (yz)^p = (zx)^q = x^2=y^2=z^2 =1\>.$$
The class of the longitude in $W$ is given by 
$((xy)^{-m}(yz)^{-k}(zx)^{-l})^2.$

Finally, Trotter identifies $W$ with a well-known triangle subgroup of isometries in the hyperbolic plane. The elements $x, y, z$ are interpreted as reflections in the three sides $\xi, \eta, \zeta$ of a hyperbolic triangle with vertex angles $\pi/p, \pi/q$ and $\pi/r$, respectively. The longitude represents a nontrivial translation along $\xi$.
Since $H, U$ are characteristic subgroups, an inversion would induce an automorphism of $W$ and  the image of the longitude would translate along $\xi$ in the opposite direction as before. This would imply that some element of $W$ reverses the direction of $\xi$. Such an element, Trotter argues, cannot exist when $p, q$ and $r$ are odd and distinct.  

In order to see how the information needed for the above proof is contained in the arc core group $AC(K)$, let $N \triangleleft AC(K)$ be the normal subgroup generated by the squares of meridional elements $g_a, a\in A(D)$. Then $AC(K)/N$ is isomorphic to $G/H$. The reason is that at as discussed above, for each crossing of $D$ the quotient class of the arc core relator is the same as the quotient class of the Wirtinger relator.  A natural isomorphism  matches classes of arc core group generators with the classes of corresponding Wirtinger generators.  As a consequence we can replace $G$ with $AC(K)$ in Trotter's proof and continue from his second step.

\subsection{Two 4-component links} \label{fourcomps}

\begin{figure} [bht]
\centering
\begin{tikzpicture} 

\draw [thick, domain=-30:295] plot ({-6+(0.8)*cos(\x)}, {(0.8)*sin(\x)});
\draw [thick, domain=150:475] plot ({-4.95+(0.8)*cos(\x)}, {(0.8)*sin(\x)});
\draw [thick, domain=-30:295] plot ({2-5+(0.8)*cos(\x)}, {(0.8)*sin(\x)});
\draw [thick, domain=150:475] plot ({2-3.95+(0.8)*cos(\x)}, {(0.8)*sin(\x)});
\node at (-6,1.1) {$s$};
\node at (-4.95,1.1) {$t$};
\node at (-3,1.1) {$u$};
\node at (0.05-2,1.1) {$v$};
\node at (-4,-1.6) {$L$};
\node at (5.5,-1.6) {$L'$};

\draw [thick, domain=-30:295] plot ({3+(0.8)*cos(\x)}, {(0.8)*sin(\x)});
\draw [thick, domain=150:295] plot ({-3.95+8+(0.8)*cos(\x)}, {(0.8)*sin(\x)});
\draw [thick, domain=-30:115]  plot ({-3.95+8+(0.8)*cos(\x)}, {(0.8)*sin(\x)});
\draw [thick, domain=150:475] plot ({-2.9+8+(0.8)*cos(\x)}, {(0.8)*sin(\x)});
\draw [thick, domain=0:360] plot ({(0.8)*cos(\x)+7}, {(0.8)*sin(\x)});
\node at (3,1.1) {$w$};
\node at (8-3.95,1.1) {$x$};
\node at (8-2.9,1.05) {$y$};
\node at (8-3.95,-1.1) {$x'$};
\node at (7,1.1) {$z$};
\end{tikzpicture}
\caption{The links $L$ and $L'$.}
\label{hhufig}
\end{figure}
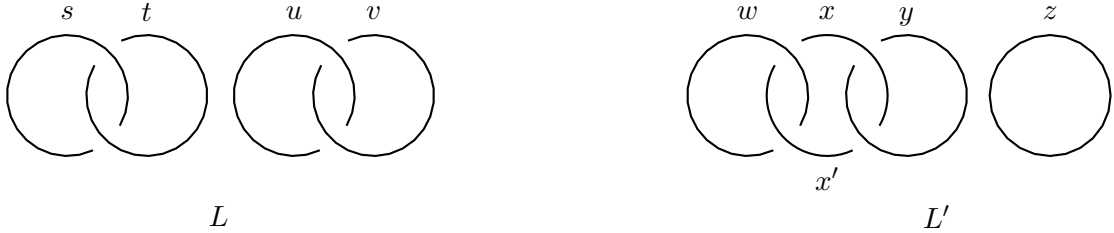

In this section we discuss the links $L=K_1 \cup K_2 \cup K_3 \cup K_4$ and $L'=K'_1 \cup K'_2 \cup K'_3 \cup K'_4$ pictured in Fig.\ \ref{hhufig}. We index their components according to the alphabetical order of the arcs. It turns out that the two links have isomorphic $AC$ groups, but no isomorphism between the groups matches the meridians to each other. In fact, the two links can be distinguished by their involutory meridional automorphisms.

In $AC(L)$ we have $g_t g_s^{-1} = g_s g_t^{-1}$, so $(g_t g_s^{-1})^2=1$. Similarly, $(g_v g_u^{-1})^2=1$. Therefore
\begin{equation} \label{iso1}
AC(L) \cong \Z *\Z_2 *\Z*\Z_2 \text{,}
\end{equation}
with the free product's factors generated by $g_t, g_s g_t^{-1}, g_v$ and $g_u g_v^{-1}$ respectively. 

In $AC(L')$ we have $g_{x'}g_w^{-1}=g_wg_{x'}^{-1}$ and $g_wg_x^{-1}=g_{x'}g_w^{-1}$; therefore $g_{x'}g_x^{-1} = g_{x'}g_w^{-1}g_wg_x^{-1}=g_wg_{x'}^{-1}g_{x'}g_w^{-1}=1$. That is, $g_x=g_{x'}$. With this equality in mind, the crossing relations of $L'$ tell us that $g_w g_x^{-1} = g_x g_w^{-1}$ and $g_y g_x ^{-1} = g_x g_y^{-1}$. Therefore $(g_w g_x^{-1})^2=1$ and $(g_y g_x^{-1})^2=1$, so
\begin{equation} \label{iso2}
AC(L') \cong \Z *\Z_2 *\Z_2*\Z \text{,}
\end{equation} 
with the free product's factors generated by $g_x, g_w g_x^{-1}, g_y g_x^{-1}$ and $g_z$ respectively. 

\begin{prop} \label{noiso}
The groups $AC(L)$ and $AC(L')$ are isomorphic, but there is no isomorphism under which $\mathcal{I}(L)$ corresponds to $\mathcal{I}(L')$.
\end{prop}

\begin{proof}
To be precise, the proposition asserts that there is no isomorphism $f:AC(L) \to AC(L')$ such that $\mathcal {I}(L') = \{f \circ \i \circ f^{-1} \mid \i \in \mathcal{I}(L) \}$.  

To verify this assertion, we develop a standard form for elements of $\mathcal{I}(L)$. Notice first that for any $b \in A(D)$,
\[
\i_{g_s}\i_{g_t}(g_b) =\i_{g_s}(g_t g_b^{-1} g_t) =g_s g_t^{-1} g_b g_t^{-1} g_s
= g_t g_s ^{-1} g_b g_s^{-1} g_t   = \i_{g_t}(g_s g_b^{-1} g_s)=\i_{g_t}\i_{g_s}(g_b) \text{,}
\]
so $\i_{g_s} \circ \i_{g_t} = \i_{g_t} \circ \i_{g_s}$. Similarly, $\i_{g_u} \circ \i_{g_v} = \i_{g_v} \circ \i_{g_u}$. Taking into account the fact that $\i_m$ maps are involutions, these commutativity relations imply that every element $\i \in \mathcal{I}(L)$ can be written as a composition $\i = \i_1 \circ \dots \circ \i_k$ for some $k \geq 1 \in \N$, where $\i_1 \in \{1, \i_{g_s}, \i_{g_t}, \i_{g_s}\i_{g_t } \}$, $\i_j \in \{\i_{g_u}, \i_{g_v}, \i_{g_u}\i_{g_v} \}$ whenever $j$ is even, and $\i_j \in \{ \i_{g_s}, \i_{g_t}, \i_{g_s}\i_{g_t} \}$ whenever $j \geq 3$ is odd. 

Notice that $\i_{g_s}(g_s)=g_s=\i_{g_t}(g_s)$, $\i_{g_s}(g_t)=g_t=\i_{g_t}(g_t)$, and neither $\i_{g_s}$ nor $\i_{g_t}$ has a fixed point outside the subgroup of $AC(L)$ generated by $g_s$ and $g_t$. Similarly, $\i_{g_u}(g_u)=g_u=\i_{g_v}(g_u)$, $\i_{g_u}(g_v)=g_v=\i_{g_v}(g_v)$, and neither $\i_{g_u}$ nor $\i_{g_v}$ has a fixed point outside the subgroup of $AC(L)$ generated by $g_u$ and $g_v$. Now, let $\i = \i_1 \circ \dots \circ \i_k$ be an arbitrary element of $\mathcal{I}(L)$, written in the standard form of the preceding paragraph. Then the fixed points of $\i$ constitute a subgroup of $AC(L)$ that is of one of these three types: (1) all of $AC(L)$ (if $k=1$ and $\i_1=1$), (2) the trivial subgroup $\{1\}$ (if $k>2$, or $k=2$ and $\i_1 \neq 1$), or (3) isomorphic to $\Z * \Z_2$ (if $k=1$ and $\i_1 \neq 1$, or $k=2$ and $\i_1=1$). 

Of course it follows that if $f:AC(L) \to AC(L')$ is an isomorphism then for an arbitrary element $\i \in \mathcal{I}(L)$, the fixed points of $f \circ \i \circ f^{-1}$ constitute a subgroup of $AC(L')$ that is all of $AC(L')$, or the trivial subgroup $\{1\}$, or isomorphic to $\Z * \Z_2$. Therefore $\i_{g_x}$, whose fixed points include $g_w, g_x$ and $g_y$ but not $g_z$, is not equal to $f \circ \i \circ f^{-1}$ for any $\i \in \mathcal{I}(D)$.\end{proof}

\begin{cor} \label{noiso1}
The groups $AC(L)$ and $AC(L')$ are isomorphic, but there is no isomorphism between them that matches meridians to meridians.
\end{cor}
\begin{proof} There are at least two different ways to show that there is no isomorphism $f:AC(L) \to AC(L')$ with $f(\M(L)) = \M(L')$. The first is that $\mathcal{I}(L)$ would necessarily correspond to $\mathcal{I}(L')$ under $f$, contradicting Proposition \ref{noiso}. The second is that $f$ would induce an isomorphism between $O(L)$ and $O(L')$. This is impossible because (\ref{iso1}) implies $O(L) \cong (\Z_2 \oplus \Z_2)*(\Z_2 \oplus \Z_2)$, while (\ref{iso2}) implies $O(L') \cong O(L' \setminus K'_4)*\Z$. \end{proof}

We take the time to prove Proposition \ref{noiso} and Corollary \ref{noiso1} to make the point that $\mathcal{I}(L)$ and $\M(L)$ can carry information even if we pay no attention to the longitudes. However, it is much easier to distinguish $L$ from $L'$ using longitudes: all four longitudes of $L$ are nontrivial, but the longitude of $K'_4$ in $AC(L')$ is trivial.

The links $L$ and $L'$ illustrate a point mentioned in the introduction: the meridians and longitudes in meridian-longitude pairs are independent. In $L$, notice that for any orientations of the components, $\la(s) = g_tg_s^{-1}$ and $\la(t) = g_sg_t^{-1}$. As $(g_tg_s^{-1})^2 = 1$, $\la(s)=\la(t)$ so the meridian-longitude pairs $(g_s, \la(s))$ and $(g_t,\la(t))$ share a longitude but have distinct meridians; indeed the meridians represent different components. Now, give $K_2'$ the counterclockwise orientation. Then $\la(x)$ and $\la(x')$ both appear in meridian-longitude pairs with $g_x=g_{x'}$; the isomorphism (\ref{iso2}) guarantees that $\la(x)=g_w g_y^{-1} = g_w g_x^{-1} \cdot 
 (g_y g_x^{-1})^{-1}= g_w g_x^{-1} \cdot 
 g_y g_x^{-1} \neq g_y g_x^{-1} \cdot 
 g_w g_x^{-1}  = g_y g_w^{-1} = \la(x')$.

 Another property of core groups illustrated by $L$ and $L'$ was mentioned in Section \ref{ab}: to find the core group of a sublink obtained by removing one link component we must know all the meridians, not just a meridian of the component that is removed. $L$ and $L'$ illustrate this property because (\ref{iso1}) and (\ref{iso2}) provide an isomorphism $AC(L) \cong AC(L')$ under which a meridian of $K_4$ corresponds to a meridian of $K'_4$; but $AC(L \setminus K_4) \not \cong AC(L' \setminus K'_4)$.

 \subsection{A family of links with homeomorphic complements} \label{homeo}

\begin{figure}[H]
\begin{center}
\includegraphics[height=2.3 in]{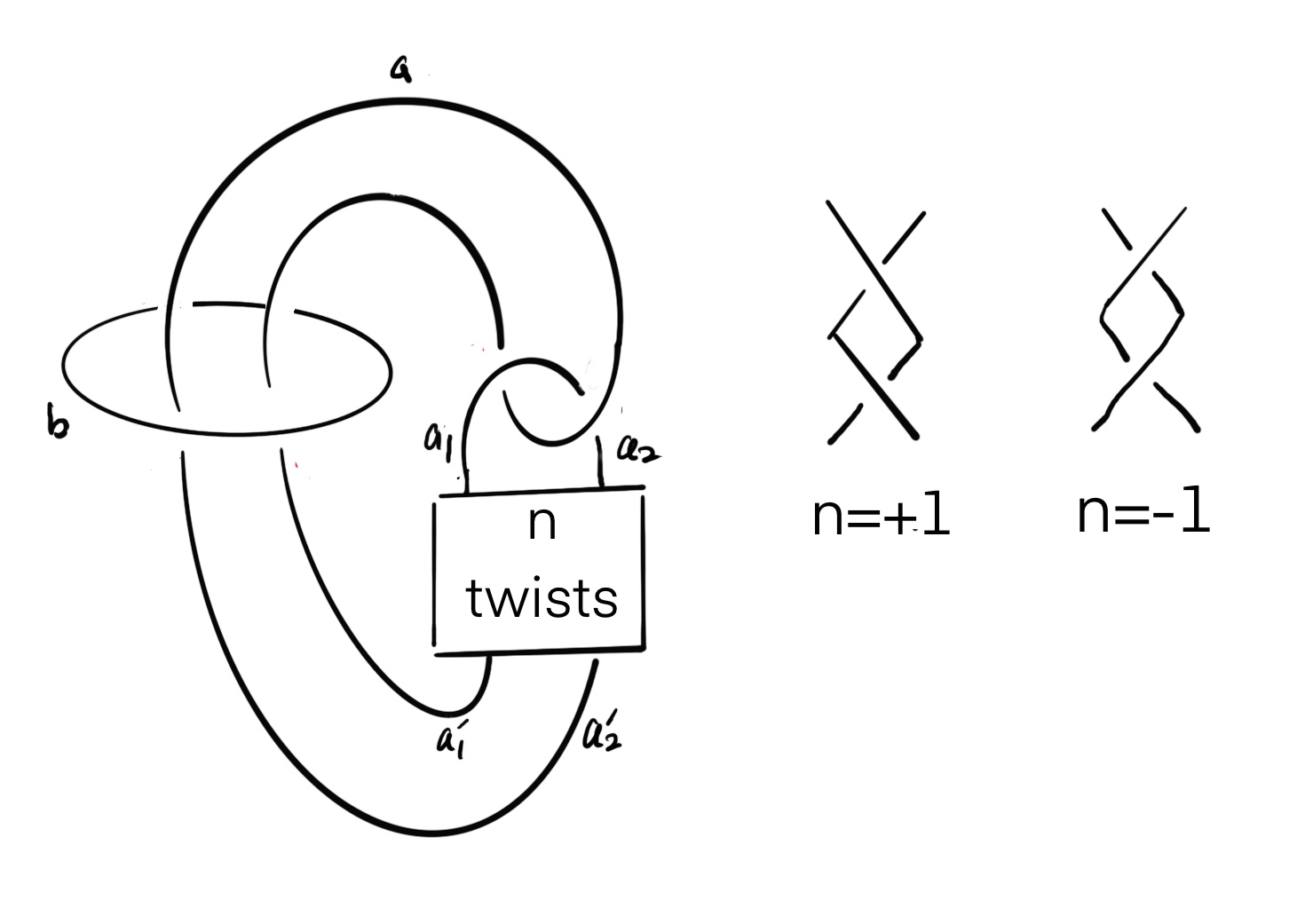}
\caption{Links $W_n$ (left); Single full twist (right)}
\label{twisted}
\end{center}
\end{figure}

Consider the infinite family of distinct 2-component links $W_n, n \in \Z$, obtained from the Whitehead link $W_0$ by giving $n$ full twists to the two strands passing through one component, as in Figure \ref{twisted}. Since every $W_n$ can also be described as the result of Dehn surgery along an unknot, the complements of the $W_n$ are all homeomorphic. In order to distinguish the links using the fundamental groups of the complements it would be necessary to introduce peripheral structures. We show that this information is already encoded in the core groups. 

\begin{prop} \label{coretwisted} $AC^{\pm}(W_n)$ is a semidirect product $\Z/8 \ltimes \Z/(|4n-1|)$ with a generator of $\Z/8$ acting on $\Z/(|4n-1|)$ by inversion. 
\end{prop}

\begin{proof} By Proposition \ref{pm}, the core group $AC(W_n)$ is the internal free product of an infinite cyclic group generated by $g_{a_2}$ and $AC^{\pm}(W_n)$ (see Figure \ref{twisted}). The Kurosh Subgroup Theorem \cite{mks} implies that the group $AC^{\pm}(W_n)$ is an invariant of the link $W_n$. Invariance can also be seen topologically: the group $AC^{\pm}(W_n)$ is isomorphic to the fundamental group of the 2-fold branched cover of the link \cite{S, W}. We compute a presentation of $AC^{\pm}(W_n)$ using the core group generators of Figure \ref{twisted} and setting $g_{a_2}$ equal to the identity. As a convenient abuse of notation, for a diagram arc labeled $x$ we use the same label $x$ to denote the image of $g_x$ after $g_{a_2}$ is modded out.

It is easy to verify that $a_1=a^2$. By induction we have $a_1' = a_1^{2n+1} = a^{4n+2}$ and $a_2'=a_1^{2n}=a^{4n}$. Then simplifying the presentation of $AC^{\pm}(W_n)$ using Tietze transformations produces the following.

$$AC^{\pm}(W_n) \cong \langle a, b ; a^2ba^2b^{-1}, a^{4n-1}b^2\rangle.$$

\begin{lemma}\label{order} The generator $b$ has order 8 in $AC^{\pm}(W_n)$. \end{lemma}

\begin{proof} (Cornelius Pillen suggested the following argument to the authors.) 
Let $m= 4n-1$. The second relator $a^{4n-1}b^2$ implies  that ({\it i}) $a^{2m} = b^{-4}$, while  the first implies that ({\it ii}) $ba^{2m}b^{-1} = a^{-2m}$. Substituting $b^{-4}$ for $a^{2m}$ in ({\it ii}), we obtain $b b^{-4}b^{-1} = b^4$. 
Hence the order of $b$ divides 8. Since the order of $b$ is equal to 8 in the abelianization of $AC^{\pm}(W_n)$,  its order is  equal to 8 in the group $AC^{\pm}(W_n)$.

\end{proof}

Continuing the proof of Proposition \ref{coretwisted}, we introduce a generator $u$ and defining relation $u=ab^{-2}$. Simplifying, we have:

$$AC^{\pm}(W_n)\cong \langle b, u ; ub^2ubub^2ub^3, (ub^2)^{4n-1} b^2\rangle.$$

Consider the short exact sequence 
$$1 \rightarrow K \rightarrow AC^{\pm}(W_n) \xrightarrow{\chi} \Z/8 \rightarrow 1,$$
where $\chi(b) =1$ and $\chi(u)=0$.  The kernel $K$ has presentation: 

\begin{equation}\label{pres}
K \cong \langle u_j ; u_ju_{j+2}u_{j+3}u_{j+5}, \quad u_ju_{j+2}\cdots u_{j+|8n-2|-2}\rangle,
\end{equation}
where the index $j$ ranges over the integers modulo 8. The relators in the second set are words of length $|4n-1|$ 
(with generators repeating).

The presentation of $K$ can be obtained from the Reidemeister-Schreier method \cite{mks}. However, it also can be found directly. That the elements $u_j$ generate is seen by considering an arbitrary product $u^{i_i}b^{j_1}\cdots u^{i_m}b^{j_m}$ in $AC^{\pm}(W_n)$ and noting that it is in the kernel of $\chi$ if and only if $\sum j_k=0\ {\rm mod}\ 8$. Lemma \ref{order} implies that such a product can be rewritten as a product of elements $b^jub^{-j}$ denoted by $u_j$. To see that the relations suffice, begin with a standard 2-complex for the presentation of $AC^{\pm}(W_n)$ and lift the 2-cells corresponding to the two relators. 
Details are left to the reader. 

The presentation (\ref{pres}) consists of two sets of relators, each consisting of 8 members. The first set can be written: 
$$u_0u_2 = u_5^{-1}u_3^{-1}$$
$$u_1u_3 = u_6^{-1}u_4^{-1}$$
$$u_2u_4 = u_7^{-1}u_5^{-1}$$
$$u_3u_5 = u_0^{-1}u_6^{-1}$$
$$u_4u_6 = u_1^{-1}u_7^{-1}$$
$$u_5u_7 = u_2^{-1}u_0^{-1}$$
$$u_6u_0 = u_3^{-1}u_1^{-1}$$
$$u_7u_1 = u_4^{-1}u_2^{-1}.$$
The reader can check that these imply that  
$$u_0u_2=u_2u_4=u_4u_6=u_6u_0 = (u_1u_3)^{-1} = (u_3u_5)^{-1} = (u_5u_7)^{-1} = (u_7u_1)^{-1}.$$
Denote $u_0u_2$ by $A$; then $u_1u_3 = A^{-1}$. 

Each relator in the second set is a product of $|4n-1|$ generators. In each all generators but the last can be paired off, with each pair equal to $A$ or $A^{-1}$. The number of pairs is $k = \frac{1}{2}(|4n-1|-1)$. The relators can then be written: 

$$A^k u_ 4$$
$$A^{-k} u_ 5$$
$$A^k u_ 6$$
$$A^{-k} u_ 7$$
$$A^k u_ 0$$
$$A^{-k} u_ 1$$
$$A^k u_ 2$$
$$A^{-k} u_ 3$$

These relators are equivalent to $u_0 = u_2 = u_4 = u_6 = A^{-k}$ and $u_1 = u_3 = u_5 = u_7 = A^k = u_0^{-1}$. 
Hence $K$ is generated by $u_0$, and the presentation (\ref{pres}) gives us:

$$K\cong \langle u_0; u_0^0,  u_0^{|4n-1|}, u_0^{-|4n-1|}\rangle  \cong \langle u_0 ; u_0^{|4n-1|}\rangle \cong \Z/(|4n-1|).$$

Moreover, $bu_0b^{-1} = u_1 = u_0^{-1}$. 
\end{proof} 

\begin{remark} The first homology groups of the 2-fold branched cyclic covers of the links $W_n$ are all isomorphic to 
$\Z /8$. However, the orders of the fundamental groups, $AC^{\pm}(W_n)$, grow linearly with respect to $|n|$. 

\end{remark}

\subsection{The Borromean rings}\label{Borr}

\begin{figure} [bth]
\centering
\begin{tikzpicture} [>=angle 90]
%left
\draw [thick, domain=-240:-90] plot ({-4+1.5*cos(\x)}, {-1.4+1.5*sin(\x)});
\draw [<-] [thick, domain=-90:-10] plot ({-4+1.5*cos(\x)}, {-1.4+1.5*sin(\x)});
\draw [thick, domain=10:100] plot ({-4+1.5*cos(\x)}, {-1.4+1.5*sin(\x)});
\draw [thick, domain=10:180] plot ({-5+1.5*cos(\x)}, {1.5*sin(\x)});
\draw [<-] [thick, domain=180:240] plot ({-5+1.5*cos(\x)}, {1.5*sin(\x)});
\draw [thick, domain=260:350] plot ({-5+1.5*cos(\x)}, {1.5*sin(\x)});
\draw [thick, domain=236:360] plot ({-3+1.5*cos(\x)}, {1.5*sin(\x)});
\draw [thick, domain=140:220] plot ({-3+1.5*cos(\x)}, {1.5*sin(\x)});
\draw [<-] [thick, domain=0:124] plot ({-3+1.5*cos(\x)}, {1.5*sin(\x)});
\node at (-6.5,-2.7) {$B$};
\node at (-4.7,-1.7) {$v$};
\node at (-6.2,.5) {$u$};
\node at (-4.7,.5) {$x$};
\node at (-1.9,.5) {$w$};
\node at (-2.6,-.3) {$z$};
\node at (-4.5,-2.5) {$y$};
\node at (-6.4,1.2) {$K_1$};
\node at (-1.6,1.2) {$K_2$};
\node at (-2.7,-2.7) {$K_3$};
%right
\draw [->] [thick, domain=-240:-90] plot ({2.5+1.5*cos(\x)}, {-1.4+1.5*sin(\x)});
\draw [thick, domain=-90:-10] plot ({2.5+1.5*cos(\x)}, {-1.4+1.5*sin(\x)});
\draw [thick, domain=10:100] plot ({2.5+1.5*cos(\x)}, {-1.4+1.5*sin(\x)});
\draw [thick, domain=10:180] plot ({1.5+1.5*cos(\x)}, {1.5*sin(\x)});
\draw [<-] [thick, domain=180:240] plot ({1.5+1.5*cos(\x)}, {1.5*sin(\x)});
\draw [thick, domain=260:350] plot ({1.5+1.5*cos(\x)}, {1.5*sin(\x)});
\draw [thick, domain=236:360] plot ({3.5+1.5*cos(\x)}, {1.5*sin(\x)});
\draw [thick, domain=140:220] plot ({3.5+1.5*cos(\x)}, {1.5*sin(\x)});
\draw [<-] [thick, domain=0:124] plot ({3.5+1.5*cos(\x)}, {1.5*sin(\x)});
\node at (5,-2.7) {$B'$};
\node at (1.8,-1.7) {$v$};
\node at (.3,.5) {$u$};
\node at (1.8,.5) {$x$};
\node at (4.6,.5) {$w$};
\node at (3.9,-.3) {$z$};
\node at (2,-2.5) {$y$};
\node at (0.1,1.2) {$K'_1$};
\node at (4.9,1.2) {$K'_2$};
\node at (1.2,-2.7) {$K'_3$};
\end{tikzpicture}
\caption{Two oriented versions of the Borromean rings.}
\label{borr}
\end{figure}
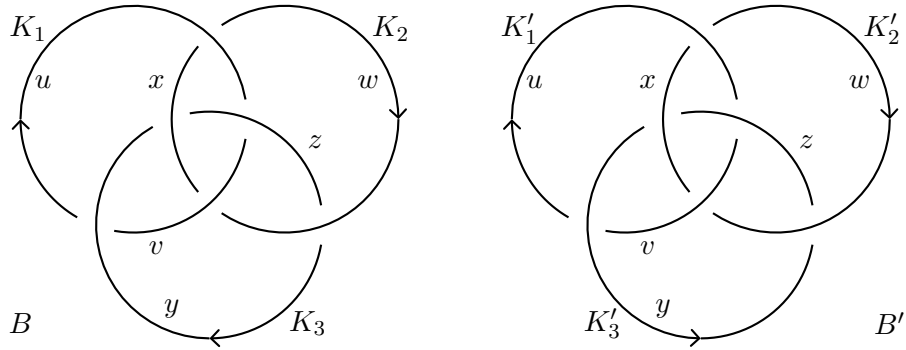

Two indexed, oriented versions of the Borromean rings are pictured in Fig.\ \ref{borr}. It is well known that the indexed, oriented link type of the Borromean rings is not changed by taking the mirror image, or by reversing the orientations of two components, so the two pictured diagrams $B$ and $B'$ represent the only two indexed, oriented versions of the link. It is not easy to verify that these two versions are actually distinct, because many link invariants do not distinguish them; for instance, their fundamental quandles are isomorphic. The two links are distinguished by Milnor's $\bar \mu$-invariants \cite{Mi}, and by the peripheral structures of their reduced Alexander modules \cite{T1}. In this section we show that they are also distinguished by the peripheral structures of their core groups.

The definition of $AC$ does not involve component orientations, so the groups $AC(B)$ and $AC(B')$ are identical. The identity map gives the impression of disrespecting the peripheral structures, because $(g_y, g_x g_w^{-1})$ is a meridian-longitude pair of $B$ while $(g_y, g_w g_x^{-1})$ is a meridian-longitude pair of $B'$. It takes a bit of work to verify this impression, i.e., to prove that $(g_y, g_w g_x^{-1})$ is not a meridian-longitude pair in $AC(B)$ and $(g_y, g_x g_w^{-1})$ is not a meridian-longitude pair in $AC(B')$. After doing this, we prove that there is no automorphism of $AC(B)$ that maps the peripheral structure of $AC(B)$ to the peripheral structure of $AC(B')$.

\begin{lemma} \label{inforder}
    The longitudes $\la(u), \la(v), \la(w), \la(x),\la(y), \la (z)$ are of infinite order in $AC(B)$.
\end{lemma}

\begin{proof}
Let $n \in \N$ be a positive integer. We use standard cycle notation for nontrivial elements of the symmetric group $S_{2n}$, and 1 for the identity. The reader can easily check that there is a homomorphism $\eta_n:AC(B) \to S_{2n}$ given by the following assignments.
\begin{align*}
    &\eta_n(g_u)=1 \\
    &\eta_n(g_v) = (1 3 5 \dots (2n-1))(2 4 6 \dots (2n)) \\
    &\eta_n(g_w) = ((2n-1) 2 (2n-3) 4 \dots 3 (2n-2) 1 (2n)) \\
    & \eta_n(g_x) = \eta_n(g_w)^{-1} = ((2n) 1 (2n-2) 3 \dots 4 (2n-3) 2 (2n-1)) \\
    &\eta_n(g_y) = \eta_n(g_z) = (1 2 3 \dots (2n))
\end{align*}

The image of $\la(w) = g_vg_u^{-1}$ is $\eta_n(g_v)$, whose order is $n$. As $n$ is arbitrary, the order of $\la(w)$ itself must be infinite. Symmetry assures us that the same is true of the other longitudes. \end{proof} 

\begin{lemma} \label{sixlong}
    There are precisely six different longitudes in $AC(B)$, namely $\la(u), \la(v), \la(w), \la(x), \allowbreak \la(y)$ and $\la (z)$.
\end{lemma}

\begin{proof}
To verify that $\la(u), \dots, \la (z)$ are six different elements of $AC(B)$, note first that according to Lemma \ref{inforder}, the two corresponding to each component are distinct: $\la (u) = g_z g_y ^{-1} = \la (v)^{-1} \neq \la(v)$, $\la (w) = g_v g_u ^{-1} = \la (x)^{-1} \neq \la(x)$, and $\la (y) = g_x g_w ^{-1} = \la (z)^{-1} \neq \la(z)$. To verify that no two of $\la(u), \dots, \la (z)$ corresponding to different components are the same, note that the kernel of the homomorphism $\eta_n:AC(B) \to S_{2n}$ used in the proof of Lemma \ref{inforder} includes $\la(u)=g_z g_y^{-1}$ and $\la(v)=g_y g_z^{-1}$, but does not include any of $\la(w), \la(x),\la(y), \la (z)$. Therefore $\la(u), \la(v) \notin \{ \la(w), \la(x),\la(y), \la (z)\}$. Symmetry assures us that also $\la(w), \la(x) \notin \{ \la(u), \la(v),\la(y), \la (z)\}$.

To complete the proof we must show that no other element of $AC(B)$ is a longitude. To verify this, observe the following.
\begin{align*}
&\i_{g_u}(\la (u)) = \i_{g_u}(g_z g_y^{-1})= g_u g_z^{-1} g_u (g_u g_y^{-1} g_u)^{-1} = g_u g_z^{-1} \cdot g_y g_u^{-1} = g_z g_v^{-1} \cdot g_v g_y^{-1} = g_z g_y^{-1} = \la(u) \\
&\i_{g_u}(\la (v)) = \i_{g_u}(g_y g_z^{-1})= g_u g_y^{-1} g_u (g_u g_z^{-1} g_u)^{-1} = g_u g_y^{-1} \cdot g_z g_u^{-1} = g_y g_v^{-1} \cdot g_v g_z^{-1} = g_y g_z^{-1} = \la(v) \\
&\i_{g_u}(\la (w)) = \i_{g_u}(g_v g_u^{-1})= g_u g_v^{-1} g_u (g_u g_u^{-1} g_u)^{-1} = g_u g_v^{-1} = \la(x) \\
&\i_{g_u}(\la (x)) = \i_{g_u}(g_u g_v^{-1})= g_u g_u^{-1} g_u (g_u g_v^{-1} g_u)^{-1} = g_v g_u^{-1} = \la(w) \\
&\i_{g_u}(\la (y)) = \i_{g_u}(g_x g_w^{-1})= g_u g_x^{-1} g_u (g_u g_w^{-1} g_u)^{-1} = g_u g_x^{-1} \cdot g_w g_u^{-1} = g_w g_u^{-1} \cdot g_u g_x^{-1} = g_w g_x^{-1} = \la(z) \\
&\i_{g_u}(\la (z)) = \i_{g_u}(g_w g_x^{-1})= g_u g_w^{-1} g_u (g_u g_x^{-1} g_u)^{-1} = g_u g_w^{-1} \cdot g_x g_u^{-1} = g_x g_u^{-1} \cdot g_u g_w^{-1} = g_x g_w^{-1} = \la(y) 
\end{align*}

That is, the action of $\i_{g_u}$ on the set $\{ \la(u), \la(v), \la(w), \la(x),\la(y), \la (z)\}$ is to fix each longitude associated with the same component as $u$, and to transpose the two longitudes corresponding to each other component. Symmetry assures us that for every $a \in \{v,w,x,y,z\}$ the action of $\i_{g_a}$ on the set $\{ \la(u), \la(v), \la(w), \la(x),\la(y), \la (z)\}$ can be described in the same way; it follows that $\i_{g_a}$ preserves the set $\{ \la(u), \dots, \la (z)\}$. Therefore no other element of $AC(B)$ is a longitude.
\end{proof}

Now, let $\a:AC(B) \to AC(B)_{ab}$ be the canonical map onto the abelianization of $AC(B)$. Then $AC(B)_{ab}$ is the abelian group with generators $\a g_u, \dots, \a g_z$ and crossing relations of the form $2 \a g_{a(c)} = \a g_{b_1(c)}+\a g_{b_2(c)}$ for each crossing $c$, where $a(c)$ is the overpassing arc at $c$ and $b_1(c), b_2(c)$ are the underpassing arcs at $c$. We use $\a g_v = 2 \a g_y - \a g_u$ to eliminate $a g_v$, $\a g_x = 2 \a g_u - \a g_w$ to eliminate $a g_x$, and $\a g_z = 2 \a g_w - \a g_y$ to eliminate $a g_z$. We conclude that $AC(B)_{ab}$ is the abelian group with generators $\a g_u, \a g_w$ and $\a g_y$, subject to the following relations.
\begin{align*}
&2 \a g_y - \a g_u = 2 (2 \a g_w - \a g_y) - \a g_u \\
&2 \a g_u - \a g_w = 2(2 \a g_y - \a g_u) - \a g_u \\
&2 \a g_w - \a g_y = 2(2 \a g_u - \a g_w) - \a g_y
\end{align*}
These relations are equivalent to the requirement that $4 \a g_y = 4 \a g_w = 4 \a g_u$. It follows that $AC(B)_{ab}$ is isomorphic to $\Z \oplus \Z_4 \oplus \Z_4$, with the three direct summands generated by $\a g_u, \a g_w - \a g_u$ and $\a g_y - \a g_u$ respectively. We identify $AC(B)_{ab}$ with the direct sum.

The images of $g_u, g_w$ and $g_y$ are the ordered triples $(1,0,0), (1,1,0)$ and $(1,0,1)$, respectively. The crossing relations of $B$ imply that twelve elements of $AC(B)_{ab} = \Z \oplus \Z_4 \oplus \Z_4$ are images of meridians. The images of meridians of $K_1$ are $\a g_u=(1,0,0), \a g_v = (1,0,2), \a(g_w g_u^{-1} g_w) = (1,2,0)$ and $\a(g_w g_v^{-1} g_w) = (1,2,2)$. The images of meridians of $K_2$ are $\a g_w = (1,1,0), \a g_x = (1,3,0), \a(g_y g_w^{-1} g_y) = (1,3,2)$ and $\a(g_y g_x^{-1} g_y) = (1,1,2)$. The images of meridians of $K_3$ are $\a g_y = (1,0,1), \a g_z = (1,2,3), \a(g_u g_y^{-1} g_u) = (1,0,3)$ and $\a(g_u g_z^{-1} g_2) = (1,2,1)$.

For each meridian $m$ of $AC(B)$, the automorphism $\i_m$ of $AC(B)$ induces an automorphism $j_m$ of $AC(B)_{ab}$. Then if $m'$ is another meridian, $j_m(\a m') = \a(\i_m(m')) = \a(m(m')^{-1}m) = 2 \a (m) - \a(m')$. The value of $2 \a(m)$ is the same for all four images of meridians associated with a single component of $B$, so there are actually only three distinct $j_m$ maps. Moreover, if $m_1$ and $m_2$ are meridians then $2( \a m_1 - \a m_2)=2(\a m_2 - \a m_1)$, so
\[
j_{m_1}j_{m_2}(m) = 2 \a m_1 - 2 \a m_2 +m = 2 \a m_2 - 2 \a m_1 +m = 
j_{m_2}j_{m_1}(m) .
\]
That is, the three $j_m$ maps commute with each other.

Each  $j_m$ map acts as the identity on the set of four images of meridians associated with the same component of $B$ as $m$, and as the product of two disjoint transpositions on the set of four images of meridians associated with some other component of $B$. For instance, $j_{g_u}$ acts as the product of transpositions $((1,1,0)(1,3,0))((1,1,2)(1,3,2))$ on the images of meridians associated with $K_2$, and the product of transpositions $((1,0,1)(1,0,3))((1,2,1)(1,2,3))$ on the images of meridians associated with $K_3$.

\begin{prop} \label{parity}
    Suppose $m$ is a meridian of $K_1$ in $AC(B)$. Then for every sequence of meridians $m_1, \dots, m_p$ such that $$\i_{m_1}\i_{m_2} \dots \i_{m_p}(m)=m \text{,}$$the number of $m_1, \dots, m_p$ associated with $K_2$ is even, and the number of $m_1, \dots, m_p$  associated with $K_3$ is even.
\end{prop}

\begin{proof}
The equality $\i_{m_1}\i_{m_2} \dots \i_{m_p}(m)=m$ implies $ j_{m_1}j_{m_2} \dots j_{m_p}(\a m)=\a m.$

As mentioned above, the $j_{m_i}$ maps all commute with each other, so we may reorganize the composition $j_{m_1}j_{m_2} \dots j_{m_p}$ in such a way that the $m_i$ corresponding to $K_1$ are listed first, the $m_i$ corresponding to $K_2$ are listed second, and the $m_i$ corresponding to $K_3$ are listed third. If $m_i$ is a meridian of $K_1$ then $j_{m_i}$ defines the identity map on the set of meridians associated with $K_1$, so the corresponding $j_{m_i}$ may be ignored. The $m_i$ of $K_2$ all define the same $j_{m_i}$ map, which is an involution. Therefore an even number of $j_{m_i}$ maps corresponding to $K_2$ cancel each other. The same applies to the meridians of $K_3$. The result is to simplify $j_{m_1}j_{m_2} \dots j_{m_p}(\a m)=\a m$ to $$ j_{g_w}^\d j_{g_y}^\e(\a m)=\a m \text{,}$$ where $\d \in \{0,1\}$ is congruent (mod 2) to the number of $m_i$ that are meridians of $K_2$, and $\e \in \{0,1\}$  is congruent (mod 2) to the number of $m_i$ that are meridians of $K_3$.  

Here are the four possible values of $j_{g_w}^\d j_{g_y}^\e(\a m)$.
\begin{align*}
    &j_{g_w}^0 j_{g_y}^0(\a m) = \a m\\
    &j_{g_w}^1 j_{g_y}^0(\a m) = 2 \a g_w - \a m = (2,2,0)-\a m\\
    &j_{g_w}^0 j_{g_y}^1(\a m) = 2 \a g_y - \a m = (2,0,2)-\a m\\
    &j_{g_w}^1 j_{g_y}^1(\a m) = 2 \a g_w - (2 \a g_y - \a_m) = (0,2,2)+ \a m.
\end{align*}
As $m$ is a meridian of $K_1$, $\a m \in \{(1,0,0), (1,2,0), (1,0,2), (1,2,2)\}$. It follows that $(2,2,0)-\a m \neq \a m$,  $(2,0,2)-\a m \neq \a m$ and $(0,2,2)+ \a m \neq \a m$. Therefore $\d = \e = 0$.\end{proof}

Of course Proposition \ref{parity} also applies to the meridians of $K_2$ and $K_3$, \emph{mutatis mutandi}.

\begin{cor} \label{onelong}
Each meridian in $AC(B)$ has only one associated longitude.
\end{cor}
\begin{proof}
Suppose $m$ is a meridian of $K_1$ that appears in meridian-longitude pairs with both $\la (u)$ and $\la (v)$. Then there is a sequence $m_1, \dots, m_p$ of meridians from the set $\{g_u,g_v,g_w,g_x,g_y,g_z\}$, such that $$\i_{m_1} \dots \i_{m_p}(m)=m \qquad \text{      and      } \qquad \i_{m_1} \dots \i_{m_p}(\la (u))=\la (v).$$ 
        
According to Proposition \ref{parity}, $\i_{m_1} \dots \i_{m_p}(m)=m$ implies that the number of $m_i$ associated with $K_2$ is even, and the number of $m_i$ associated with $K_3$ is also even.

In the proof of Lemma \ref{sixlong} we observed that for each meridian $m_i \in \{g_u,g_v,g_w,g_x,g_y,g_z\}$, the action of $\i_{m_i}$ on the set of six longitudes in $AC(B)$ may be described as follows: $\i_{m_i}$ fixes each of the two longitudes of the component associated with $m$, and $\i_{m_i}$ transposes the two longitudes associated with each of the other components. Therefore $\i_{m_1} \dots \i_{m_p}(\la (u))=\la (v)$ implies that the total number of $m_i$ associated with $K_2$ or $K_3$ is odd. This is a contradiction.
        
Symmetry assures us that the same kind of argument holds if $m$ is a meridian of $K_2$ or $K_3$. \end{proof}

In light of Corollary \ref{onelong}, when $m$ is a meridian of $AC(B)$ we can use the notation $\la(m)$ for the unique longitude of $m$ in $AC(B)$.

Lemma \ref{inforder} and Corollary \ref{onelong} tell us that the identity map of $AC(B)=AC(B')$ is not compatible with the peripheral structures of $AC(B)$ and $AC(B')$: the identity map sends $g_y \in AC(B)$ to $g_y \in AC(B')$, but it does not map the unique longitude of $g_y$ in $AC(B)$ to the unique longitude of $g_y$ in $AC(B')$. 

In the rest of this section we show that in fact, there is no automorphism of $AC(B) =AC(B')$ that maps the peripheral structure of $AC(B)$ to the peripheral structure of $AC(B')$. 

\begin{definition} \label{special}
An automorphism of $AC(B)$ is \emph{special} if it maps the peripheral structure of $AC(B)$ to the peripheral structure of $AC(B')$. \end{definition}

That is, a special automorphism $\s$ of $AC(B)$ has these properties: if $m$ is a meridian associated with $K_i$ then $\s(m)$ is also a meridian associated with $K_i$; if $m$ is a meridian associated with $K_1$ or $K_2$ then $\s(\la(m)) = \la(\s(m))$; and if $m$ is a meridian associated with $K_3$ then $\s(\la(m)) = \la(\s(m))^{-1}$.

Some algebraic properties of the longitudes will be handy in proving that special automorphisms do not exist. 

\begin{lemma}
The longitudes commute with each other in $AC(B)$.
\end{lemma}
\begin{proof}
The lemma follows from the crossing relations in $AC(B)$. For instance:
    \[
\la(u)\la(x) = g_z g_y^{-1} g_u g_v^{-1} = g_z g_v^{-1} g_y g_v^{-1} = g_u g_z^{-1} g_u g_y^{-1} = g_u g_v^{-1} g_z g_y^{-1} = \la(x) \la(u). \qedhere
    \]
\end{proof}

\begin{lemma} \label{canon}
If $m$ is a meridian of $K_3$ in $AC(B)$ then there are unique integers $r,s \in \Z$ such that $m=\la (u)^r \la (w)^s g_y$.
\end{lemma}
\begin{proof}
We begin by proving that such $r$ and $s$ exist. 

There are meridians $m_1, \dots, m_k \in AC(B)$ such that $\i_{m_1} \dots \i_{m_k}(g_y)=m$. As the group $\mathcal{I}(D)$ is generated by $\i_{g_u}, \dots, \i_{g_z}$, we may presume that for each $i \in \{1, \dots, k\}$, there is an arc $a_i \in 
 \{u,v,w,x,y,z\}$ such that $m_i = g_{a_i}$. 

The crossing relations of $B$ yield relations among these $\i_{m_i}$ maps; for instance $g_u g_w^{-1} = g_x g_u^{-1}$ implies that $\i_{g_u}\i_{g_w}=\i_{g_x}\i_{g_u}$, because for any meridian $n$, $$\i_{g_u}\i_{g_w}(n) = \i_{g_u}(g_w n^{-1} g_w) =  g_u g_w^{-1} n g_w^{-1} g_u = g_x g_u^{-1} n g_u^{-1} g_x = \i_{g_x}\i_{g_u}(n).$$ Using these relations, we can rewrite the composition $\i_{m_1} \dots \i_{m_k}$ so that the first few $m_i$ are all either $g_y$ or $g_z$, the next few $m_i$ are all either $g_u$ or $g_v$, and the last few $m_i$ are all either $g_w$ or $g_x$. As $\i_{g_w}(g_y)=g_z=\i_{g_x}(g_y)$, the effect of these last few $\i_{m_i}$ is simply to toggle back and forth between $g_y$ and $g_z$. We conclude that $$m=\i_{m_1} \dots \i_{m_{p+q}}(m')$$ where $m', m_1, \dots, m_p \in \{g_y,g_z \}$ and $m_{p+1}, \dots, m_{p+q}\in \{g_u,g_v \}$. As each $\i_{m_i}$ map is an involution, $\i_{m_i} \circ \i_{m_i}$ is always the identity map; hence we may presume that $m_i \neq m_{i+1} \thickspace \forall i$.

Now, let $g_1 ,g_2,g'_1, g'_2\in AC(B)$ be the elements
\begin{align*}
&g_1 = m_1 m_2 ^{-1} \dots m_p^{(-1)^{p+1}} \\
&g'_1 = m_p^{(-1)^{p+1}} m_{p-1}^{(-1)^{p}} \dots m_1\\ 
&g_2 = m_{p+1}^{(-1)^{p+2}} \dots m_{p+q}^{(-1)^{p+q+1}}\\
&g'_2 = m_{p+q}^{(-1)^{p+q+1}}m_{p+q-1}^{(-1)^{p+q}} \dots m_{p+1}^{(-1)^{p+2}}.
\end{align*}
Then $g_1$ is a doubly alternating product of $g_y$ and $g_z$; we say ``doubly alternating'' because $g_y$ and $g_z$ alternate, and the exponents $+1$ and $-1$ alternate. Also, $g'_1$ is the product of the same factors as $g_1$, in the opposite order. Similarly, $g_2$ is a doubly alternating product of $g_u$ and $g_v$, and $g'_2$ is the product of the same factors in the opposite order. Then the equality $m=\i_{m_1} \dots \i_{m_{p+q}}(m')$ is equivalent to the equality $$m=g_1g_2 (m')^{(-1)^{p+q+2}}g'_2g'_1.$$

Suppose for the moment that $p$ is even. Then $g_1$ is the concatenation of $p/2$ occurrences of $$(g_y \text{  or  }  g_z) \cdot (g_z \text{  or  } g_y)^{-1} \text{,}$$ so $g_1 = \la(u) ^ {\pm p/2}$. Recall that $m'$ is $g_y$ or $g_z$, and $g'_2$ is a product of factors ($g_u$ or $g_v$)$^{\pm 1}$. $AC(B)$ has crossing relations like $g_y g_u^{-1}=g_v g_y^{-1}$ and $g_z^{-1} g_v = g_u^{-1} g_z$, which allow us to push the factor $(m')^{(-1)^{p+q+2}}$ through $g'_2$ one factor at a time, resulting in an equality $$m=g_1 g_ 2 g''_2 (m')^{(-1)^{p+2}} g'_1 \text{,}$$
where $g''_2$ is obtained from $g'_2$ by changing each occurrence of ($g_u$ or $g_v$)$^{\pm 1}$ to ($g_v$ or $g_u$)$^{\mp 1}$. The product $g_2 g''_2$ is the concatenation of $q$ occurrences of $$(g_u \text{  or  }  g_v) \cdot (g_v \text{  or  } g_u)^{-1} \text{,}$$
so $g_2 g''_2 = \la(w)^{\pm q}$. The product $(m')^{(-1)^{p+2}} g'_1$ is the concatenation of $p/2$ occurrences of $$(g_y \text{  or  }  g_z) \cdot (g_z \text{  or  } g_y)^{-1} \text{,}$$followed by a single occurrence of $g_y$ or $g_z$. If the final factor is $g_z$, we may replace it with $g_z g_y^{-1} g_y = \la(u) g_y$. Also, some of the $p/2$ occurrences may cancel. All in all, we have $$(m')^{(-1)^{p+2}} g'_1 = \la(u) ^ n \cdot g_y$$ for some integer $n$, and hence
\[
m=g_1 g_ 2 g''_2 (m')^{(-1)^{p+2}} g'_1 = \la(u) ^ {\pm p/2} \cdot \la(w)^{\pm q} \cdot \la(u) ^ n \cdot g_y= \la(u) ^ {n \pm p/2} \cdot \la(w)^{\pm q}  \cdot g_y.
\]

If $p$ is odd, then after obtaining $$m=g_1 g_ 2 g''_2 (m')^{(-1)^{p+2}} g'_1$$ as above, we can take the last term in the product $g_1$ and push it through $g_2 g''_2$ using crossing relations. The result is an equality $$m=\la (u) ^{\pm (p-1)/2} \la(w) ^{\pm q} \la(u)^n g_y$$ for some $n \in \Z$.

It remains to verify uniqueness. Suppose $p,q,r,s \in \Z$ and $\la (u)^p\la (w)^q g_y = \la (u)^r \la (w)^s g_y$. Then $\la(u) ^ {p-r} = \la(w) ^ {s-q}$. Let $n \in \N$, and consider the homomorphism $\eta_n:AC(B) \to S_{2n}$ of Lemma \ref{inforder}. The kernel of $\eta_n$ contains $\la(u)$, so it also contains $\la(w) ^ {s-q}$. As $n$ is arbitrary and $\eta_n(\la(w))$ is of order $n$ in $S_{2n}$, it follows that $s-q=0$. Then $\la(u) ^ {p-r} = \la(w) ^ {s-q}=1$ implies $p-r=0$, because $\la (u)$ is of infinite order in $AC(B)$. \end{proof}

There are analogues of Lemma \ref{canon} for the other components of $B$: every meridian of $K_1$ is equal to a unique product $\la (w)^r \la (y)^s g_u$, and every meridian of $K_2$ is equal to a unique product $\la (u)^r \la (y)^s g_w$.
\begin{lemma} \label{blemma1}
Let $\s$ be a special automorphism of $AC(B)$ with $\s(g_u)=g_u$ and $\s(g_v)=g_v$. Then there is an $r \in \Z$ with $\s(g_y)=(g_z g_y ^{-1})^r g_y$ and $\s(g_z)=(g_z g_y ^{-1})^r g_z$. If $r$ is even then $\s (\la(g_y)) = \la(g_y)^{-1}$, and if $r$ is odd then $\s(\la (g_y))=\la (g_y)$.
\end{lemma}
\begin{proof}
The hypothesis that $\s$ is special implies that $\s(g_y)$ is a meridian of $K_3$ in $AC(B)$. Lemma \ref{canon} tells us that $\s(g_y) = \la (u)^r \la (w)^s g_y= (g_z g_y ^{-1})^r (g_v  g_u ^{-1})^s g_y$ for some $r,s \in \Z$. Therefore we have
\[
g_v= \s(g_v) = \s(g_y g_u^{-1} g_y) = (g_z g_y ^{-1})^r (g_v  g_u ^{-1})^s g_y g_u^{-1} (g_z g_y ^{-1})^r (g_v  g_u ^{-1})^s g_y.
\]
If $n \in \N$ then the image of this equality under the homomorphism $\eta_n$ of Lemma \ref{inforder} is
\[
\eta_n(g_v) = 1^r \eta_n(g_v)^s \eta_n(g_y) 1^{-1} 1^r \eta_n(g_v)^s \eta_n(g_y).
\]
As $\eta_n(g_v) = \eta_n(g_y) ^2$, we conclude that $\eta_n(g_y) ^2 = \eta_n(g_y) ^{4s+2}$, and hence $1=\eta_n(g_y) ^{4s}$. As $n$ is arbitrary and the order of $\eta_n(g_y)$ is $2n$, it follows that $4s=0$. Therefore $\s(g_y)=(g_z g_y ^{-1})^r g_y$.

As $\s$ is a special automorphism with $\s(g_u)=g_u$, it must be that $g_z g_y^{-1}=\la (g_u) = \s(\la (g_u)) = \s(g_z g_y^{-1})$. Therefore $$\s(g_z) = \s(g_z g_y^{-1}) \s(g_y) = g_z g_y^{-1}(g_z g_y ^{-1})^r g_y =(g_z g_y ^{-1})^r g_z g_y^{-1} g_y = (g_z g_y ^{-1})^r g_z.$$

Now, suppose $r$ is even. Then 
\[
\s(g_y) = (g_z g_y ^{-1})^r g_y =(g_z g_y ^{-1})^{r-1} g_z =  (g_z g_y ^{-1})^{\sfrac{(r-2)}{2}} g_z g_y^{-1} g_z (g_y ^{-1} g_z)^{\sfrac{(r-2)}{2}}=(\i_{g_z} \i_{g_y})^{\sfrac{(r-2)}{2}} \i_{g_z}(g_y) \text{,}
\]
the image of $g_y$ under a composition of $\i_{g_y}$ and $\i_{g_z}$ maps. As observed in the proof of Lemma \ref{sixlong}, the $\i_{g_y}$ and $\i_{g_z}$ maps fix $\la(g_y)$, so $\la(\s(g_y)) = \la(g_y)$. As $\s$ is a special automorphism and $y$ is an arc of $K_3$, $\la(\s(g_y)) \neq \s(\la(g_y))$. It follows that $\s(\la(g_y))=\la(g_y)^{-1}$.

If $r$ is odd, instead, then
\[
\s(g_y) = (g_z g_y ^{-1})^r g_y =(g_z g_y ^{-1})^{r-1} g_z =  (g_z g_y ^{-1})^{\sfrac{(r-1)}{2}} g_z (g_y ^{-1} g_z)^{\sfrac{(r-1)}{2}}=(\i_{g_z} \i_{g_y})^{\sfrac{(r-1)}{2}} (g_z) \text{,}
\]
the image of $g_z$ under a composition of $\i_{g_y}$ and $\i_{g_z}$ maps. Therefore $\la(\s(g_y)) = \la(g_z)$ and $\s(\la(g_y))=\la(g_z)^{-1}$.
\end{proof}

\begin{prop} \label{almostdone}
There is no special automorphism $\s$ of $AC(B)$ with $\s(g_u)=g_u$ and $\s(g_v)=g_v$.
\end{prop}
\begin{proof}
Suppose $\s$ is such an automorphism, with $\s(g_y)=(g_z g_y ^{-1})^r g_y$ and $\s(g_z)=(g_z g_y ^{-1})^r g_z$. Then $\s(g_w)$ is a meridian of $K_2$, so 
\[
\s(g_w) = \la (u)^p \la (y)^q g_w = (g_z g_y ^{-1})^p (g_x g_w ^{-1})^q g_w
\]
for some $p,q \in \Z$. Using crossing relations, we have
\begin{align*}
(g_z g_y ^{-1})^r g_y&= \s(g_y) = \s(g_w g_z ^{-1} g_w) \\
&= (g_z g_y ^{-1})^p (g_x g_w ^{-1})^q g_w g_z^{-1} (g_z g_y ^{-1})^{-r} (g_z g_y ^{-1})^p (g_x g_w ^{-1})^q g_w\\
&= (g_z g_y ^{-1})^p (g_x g_w ^{-1})^q g_w (g_y ^{-1} g_z)^{p-r} g_z ^{-1} (g_x g_w ^{-1})^q g_w\\
&=(g_z g_y ^{-1})^p (g_x g_w ^{-1})^q (g_z g_y ^{-1})^{p-r}  g_w (g_x ^{-1} g_w )^q g_z^{-1} g_w\\
&=(g_z g_y ^{-1})^p (g_x g_w ^{-1})^q (g_z g_y ^{-1})^{p-r}  (g_w g_x ^{-1} ) ^q g_w g_z^{-1} g_w\\
&= (g_z g_y ^{-1})^p (g_x g_w ^{-1})^{q-q} (g_z g_y ^{-1})^{p-r}  g_y \\
&=(g_z g_y ^{-1})^{2p-r}  g_y \text{,}
\end{align*}
so $(g_z g_y ^{-1})^{2p}=(g_z g_y ^{-1})^{2r}$. As $g_z g_y ^{-1} = \la(u)$ is of infinite order, it follows that $p=r$.

Also, we have
\begin{align*}
\s(g_x) &= \s(g_u g_w^{-1} g_u) =g_u g_w^{-1}  (g_x g_w ^{-1})^{-q} (g_z g_y ^{-1})^{-p} g_u\\
&= g_u (g_w^{-1} g_x )^{-q} g_w ^{-1}(g_z g_y ^{-1})^{-p} g_u\\
&=  (g_x g_w^{-1})^{-q} g_u (g_y ^{-1} g_z )^{-p} g_w ^{-1} g_u\\
&=  (g_x g_w^{-1})^{-q} (g_y g_z ^{-1})^{-p} g_u g_w ^{-1} g_u\\
&=  (g_x g_w^{-1})^{-q} (g_z g_y^{-1})^{p} g_x \text{,}
\end{align*}
and therefore
\[
\s(\la(y)) = \s(g_x g_w^{-1}) = \s(g_x) \s( g_w)^{-1}= (g_x g_w^{-1})^{1-2q}  (g_z g_y ^{-1})^{p-p} = (g_x g_w^{-1})^{1-2q} = \la(y)^{1-2q}.
\]

As $\s$ is a special automorphism, $\s(\la(y)) = \la(y)^{\pm 1}$; $\la(y)$ is of infinite order, so it must be that $q \in \{0,1\}$. In each case we find a contradiction.

If $q=1$ then $\s(\la(y)) = \la(y)^{-1}$, so Lemma \ref{blemma1} tells us that $r$ is even. We know $p=r$, so we have 
\begin{align*}
\s(g_w)&=(g_z g_y ^{-1})^p (g_x g_w ^{-1})^q g_w\\
&=(g_z g_y ^{-1})^{r} (g_x g_w ^{-1})^1 g_w\\
&=(g_z g_y ^{-1})^{r/2} (g_z g_y ^{-1})^{r/2} g_x\\
&=(g_z g_y ^{-1})^{r/2} g_x ( g_y ^{-1} g_z)^{r/2}\\
&=(\i_{g_z} \i_{g_y})^{r/2}(g_x) \text{,}
\end{align*}
the image of $g_x$ under the composition of an even number of $\i_{g_y}$ and $\i_{g_z}$ maps. According to the proof of Lemma \ref{sixlong}, it follows that $\la(\s(g_w)) = \la(g_x)$. As $\s$ is a special automorphism, it should have $\s(\la(g_w)) = \la(\s(g_w))$. But
\[
\s(\la(g_w)) = \s(g_v g_u^{-1})= \s(g_v)\s( g_u)^{-1}= g_v g_u^{-1} \neq g_u g_v^{-1} = \la(g_x).
\]

If $q=0$ a similar argument applies: $\s(\la(y)) = \la(y)$, $r$ is odd, and
\begin{align*}
\s(g_w)&=(g_z g_y ^{-1})^p (g_x g_w ^{-1})^q g_w\\
&=(g_z g_y ^{-1})^{r} (g_x g_w ^{-1})^0 g_w\\
&=(g_z g_y ^{-1})^{\sfrac{(r-1)}{2} }(g_z g_y ^{-1}) (g_z g_y ^{-1})^{\sfrac{(r-1)}{2} } g_w\\
&=(g_z g_y ^{-1})^{\sfrac{(r-1)}{2} }g_z g_y ^{-1} g_w (g_y ^{-1} g_z )^{\sfrac{(r-1)}{2} } \\
&=(g_z g_y ^{-1})^{\sfrac{(r-1)}{2} }g_z g_w ^{-1} g_z (g_y ^{-1} g_z )^{\sfrac{(r-1)}{2} } \\
&=(\i_{g_z} \i_{g_y})^{\sfrac{(r-1)}{2}} \i_{g_z}(g_w) \text{,}
\end{align*}
the image of $g_w$ under the composition of an odd number of $\i_{g_y}$ and $\i_{g_z}$ maps. The proof of Lemma \ref{sixlong} tells us that $\la(\s(g_w)) = \la(g_w)^{-1}$, but
\[
\s(\la(g_w)) = \s(g_v g_u^{-1})= \s(g_v)\s( g_u)^{-1}= g_v g_u^{-1} \neq (g_v g_u^{-1})^{-1} = \la(g_w)^{-1}.
\]
As a special automorphism, though, $\s$ should have $\s(\la(g_w)) = \la(\s(g_w))$.
\end{proof}

\begin{theorem} \label{finally}
$AC(B)$ does not have a special automorphism.
\end{theorem}
\begin{proof}
Let $\s$ be a special automorphism of $AC(B)$. Then $\s(g_u)$ is a meridian of $K_1$, so there are meridians $m_1, \dots, m_k \in AC(B)$ with $\i_{m_1} \dots \i_{m_k}(g_u)=\s(g_u)$. The composition $\s' = 
\i_{m_k} \circ \dots \circ \i_{m_1} \circ \s$ is then a special automorphism of $AC(B)$ with $\s'(g_u) = g_u$. 

The fact that $\s'$ is a special automorphism implies that $\s'(\la (g_y)) = \la(g_y)^{\pm 1}$. If $\s'(\la (g_y)) = \la(g_y)$ then
\[
\s'(g_v)g_u^{-1} = \s'(g_vg_u^{-1} )=g_v g_u^{-1}\text{,}
\]
so $\s'(g_v) =g_v$. This contradicts Proposition \ref{almostdone}, so it must be that $\s'(\la (g_y)) = \la(g_y)^{- 1}$. Then

\[
\s'(g_v)g_u^{-1} = \s'(g_vg_u^{-1} )=(g_v g_u^{-1})^{-1} = g_u g_v^{-1} \text{,}
\]
so $\s'(g_v) = g_u g_v^{-1} g_u$. Now let $\s '' = \i_{g_u} \circ \s'$. Then $\s''$ is a special automorphism with $\s''(g_u)=\i_{g_u} (\s'(g_u)) = \i_{g_u}(g_u)=g_u$ and $\s''(g_v)=\i_{g_u} (\s'(g_v)) = \i_{g_u}(g_u g_v^{-1} g_u)=g_u (g_u g_v^{-1} g_u)^{-1} g_u = g_v$, again contradicting Proposition \ref{almostdone}. \end{proof}

Theorem \ref{finally} is equivalent to the assertion that there is no isomorphism $AC(B) \cong AC(B')$ that matches the two peripheral structures. That is, the two versions of the Borromean rings in Fig. \ref{borr} are distinguished by the peripheral structures of their core groups.

\section*{ORCID} 

Daniel S. Silver https://orcid.org/0000-0003-1506-3525

\noindent Lorenzo Traldi https://orcid.org/0000-0003-1097-2818

%Susan G. Williams https://orcid.org/0009-0009-5987-3168

\bigskip

\ni Department of Mathematics and Statistics,\\
\ni University of South Alabama\\ Mobile, AL 36688 USA\\
\ni Email: silver@southalabama.edu

\bigskip

\ni Department of Mathematics,\\
\ni Lafayette College\\Easton PA 18042\\
\ni Email: traldil@lafayette.edu\\ 

\end{document}